\documentclass[reqno,12pt]{amsart}
\usepackage{amssymb,amsmath,amsthm,amsfonts,color,breqn}  
\usepackage[left=2.7cm,right=2.7cm,top=2.5cm,bottom=2.7cm]{geometry}
\usepackage{mathtools} 
\usepackage[abs]{overpic} 
\usepackage{xcolor}
 \usepackage{enumerate}
\usepackage{comment}

\newtheorem{theorem}{Theorem}[section]
\newtheorem{lemma}[theorem]{Lemma}

\newcommand{\Rz}{\mathbb{R}}
\newcommand{\Nz}{\mathbb{N}}

\newcommand{\epsi}{\varepsilon}
\newcommand{\eps}{\varepsilon}

\newcommand{\RRR}{\color{red}}
\newcommand{\EEE}{\color{black}}

\renewcommand{\lvert}{\|}
\renewcommand{\rvert}{\|}

\newcommand{\mathsout}[1]
{\bgroup\mathchoice
	{\sbox0{$\displaystyle{#1}$}%
		\usebox0\hspace{-\wd0}%
		\rule[0.5\ht0-0.5\dp0-.5pt]{\wd0}{0.1pt}}%
	{\sbox0{$\textstyle{#1}$}%
		\usebox0\hspace{-\wd0}%
		\rule[0.5\ht0-0.5\dp0-.5pt]{\wd0}{0.1pt}}%
	{\sbox0{$\scriptstyle{#1}$}%
		\usebox0\hspace{-\wd0}%
		\rule[0.5\ht0-0.5\dp0-.5pt]{\wd0}{0.1pt}}%
	{\sbox0{$\scriptscriptstyle{#1}$}%
		\usebox0\hspace{-\wd0}%
		\rule[0.05\ht0-0.5\dp0-.5pt]{\wd0}{0.1pt}}%
	\egroup}

\def\eps{\varepsilon}

\renewcommand{\d}{{\mathrm d}}

\def\to{\rightarrow}

\def\wstarto{\stackrel{*}{\rightharpoonup}}
\def\weakstar{\wstarto}
\def\embed{\hookrightarrow}

\makeatletter
\@namedef{subjclassname@2020}{\textup{2020} Mathematics Subject Classification}
\makeatother
 
 \numberwithin{equation}{section}
 
\begin{document}
\title[WED Approach to Semilinear Gradient Flows]
{Weighted Energy-Dissipation Approach to Semilinear Gradient Flows with State-Dependent Dissipation}

\author{Goro Akagi}
\address[Goro Akagi]{Mathematical Institute and Graduate School of Science, Tohoku University, Aoba, Sendai 980-8578, Japan.}
\email{goro.akagi@tohoku.ac.jp}
\urladdr{http://www.math.tohoku.ac.jp/$\sim$akagi}

\author{Ulisse Stefanelli}
\address[Ulisse Stefanelli]{Faculty of Mathematics,  University of Vienna, Oskar-Morgenstern-Platz 1, A-1090 Vienna, Austria,
Vienna Research Platform on Accelerating Photoreaction Discovery,
University of Vienna, W\"ahringer Str. 17, 1090 Vienna, Austria, $\&$ Istituto di Matematica Applicata e Tecnologie Informatiche ``E. Magenes'' - CNR, v. Ferrata 1, I-27100 Pavia, Italy.}
\email{ulisse.stefanelli@univie.ac.at}
\urladdr{http://www.mat.univie.ac.at/$\sim$stefanelli}

\author{Riccardo Voso}
\address[Riccardo Voso]{Faculty of Mathematics, University of Vienna, 
	Oskar-Morgenstern-Platz 1, A-1090 Vienna, Austria,  $\&$
        Vienna School of Mathematics, University of Vienna, 
	Oskar-Morgenstern-Platz 1, A-1090 Vienna, Austria}
\email{riccardo.voso@univie.ac.at}



\begin{abstract} 
We investigate the Weighted Energy-Dissipation variational
approach to semilinear gradient flows with state-dependent
dissipation. A family of parameter-dependent functionals defined over
entire trajectories is introduced and proved to admit
 global  minimizers. These  global  minimizers correspond to solutions  of   elliptic-in-time
regularizations of the limiting causal problem. By
passing to the limit in the parameter
we prove that such  global  minimizers converge, up to subsequences,
to a solution of the gradient flow.
\end{abstract}

\maketitle

\section{Introduction}
We are interested in a global-in-time variational approach
to 
abstract gradient flows of the form
\begin{equation}
  \label{eq:00}
{\d}_{2} \psi (u ,\dot{u} ) + Au(t) + \partial \phi^2(u )\ni
0\quad\text{in}\  H,\ \text{a.e. in}\  (0,T).
\end{equation}
Here, $u:[0,T]\to H$ is a trajectory from the finite time interval
$[0,T]$ to the Hilbert space
$H$ and the dot denotes time differentiation. In relation
\eqref{eq:00}, the smooth functional
$(u,v)\in H\times H \mapsto \psi(u,v)\in [0,\infty)$ is
assumed to be convex in the {\it rate}
$v$, the symbol ${\d}_2$ indicates the Fr\'echet differential with
respect to the second variable, and the monotone map $v \mapsto {\d}_2\psi(u,v)$ is
assumed to be nondegenerate and  linearly bounded,  uniformly for $u$. The linear, continuous, and self-adjoint operator $A: X \to X^*$
maps a second Hilbert space $X$, compact and dense in $H$, to its
dual $X^*$. The functional $\phi^2: H \to [0,\infty]$ is convex and
$\partial \phi^2$ indicates its subdifferential. The doubly
nonlinear inclusion
\eqref{eq:00} hence corresponds to the gradient flow of the {\it energy}
$\phi(u)=  \langle Au,u\rangle_X/2+ \phi^2(u)$, where $\langle\cdot,\cdot \rangle_X$ denotes the duality pairing between $X^*$ and $X$, with respect to the {\it state-dependent dissipation} $\psi$. Note in particular that the
nonlinearity $A+ \partial \phi^2$  is tailored to the variational
formulation of   a semilinear elliptic
operator and that we implicitly ask $Au\in H$.  Relation
\eqref{eq:00} hence corresponds to a {\it strong} formulation of a
semilinear parabolic problem with state-dependent dynamics, see
Section \ref{sec:application}.  

This paper is concerned with the analysis of the {\it Weighted
Energy-Dissipation} (WED) approach to the gradient flow
\eqref{eq:00}.  This  is a global-in-time variational approach, hinging
on the minimization of the $\eps$-dependent $(\eps>0)$ functionals 
$W^{\eps} : L^2(0,T;H) \rightarrow [0,\infty]$ given by
\begin{equation}
\label{WEDfunctional}
W^{\eps}(u) =
\begin{cases}
\displaystyle\int_0^{T} e^{-t/\eps} \Big( \eps\psi(u(t),\dot{u}(t)) +
\phi(u(t))  \Big) \,\d t  & \text{if}\  u\in K(u_0),\\[5mm]
\infty  & \text{else},
\end{cases}
\end{equation}
where the convex domain $K(u_0)$ depends on the initial condition $u_0
\in H$ and
is defined as
\begin{align}\label{K}
K(u_0) &= \{ u \in L^2(0,T;H) \; : \; u(0) = u_0 \; \text{and} \;
         \phi(u), \ \psi(u,\dot{u}) \in L^1(0,T)       \}.
\end{align}

The functionals $W^\eps$ feature an $\eps$-dependent exponentially
decaying {\it weight} multiplying the (weighted) sum of the {\it
  energy} $\phi$ and the
{\it dissipation} $\psi$. As such, they are usually referred to as {Weighted
Energy-Dissipation} functionals. Note that the parameter $\eps>0$
plays the role of a relaxation time, which justifies the factor
$\eps$ in front of the dissipation~$\psi$.

By means of the functionals $W^{\eps}$ one can approximate the
differential problem \eqref{eq:00} on purely variational terms. The
focus of this paper is on providing a rigorous analysis of such an approximation. Our
main result (Theorem \ref{theorem}) proves that

\begin{enumerate} 
\item[(i)] for all $\eps>0$,  there exists  a global minimizer $u^\eps$ of
  $W^\eps$,
  \item[(ii)]  $u^\eps$ solves the Euler--Lagrange
    equation strongly, 
    \item[(iii)]  up to subsequences,  $u^\eps$ converge to a
      solution of \eqref{eq:00}  as $\eps\to 0$. 
\end{enumerate}

In fact, the minimization of $W^\eps$ corresponds  to
an elliptic-in-time regularization of the gradient flow
\eqref{eq:00}.  This   is revealed by  computing the
Euler--Lagrange problem corresponding to the minimization of $W^\eps$,
namely,
\begin{align}
  &-\eps \frac{{\rm d}}{{\rm d}t} ({\d}_2\psi(u,\dot u)) + \eps {\d}_1\psi(u,\dot u) + {\d}_2
  \psi(u,\dot u) + Au + \partial\phi^2 (u)\ni 0\quad\text{in}\  H ,\label{eulerolagrangeepsilon0}\\
  &u(0)=u_0, \quad \eps {\d}_2 \psi(u(T),\dot u(T)) = 0, \label{eulerolagrangeepsilon40}
\end{align}
to be solved almost everywhere in $(0,T)$. 
Here, ${\d}_1$ denotes the Fr\'echet differential with respect to the
first variable and a natural Neumann condition is imposed at the final time $T$.
 Due to ellipticity in time, the  global  minimizers $u^\eps$ are more
regular  than the
solution to the
limiting problem \eqref{eq:00}.  On the other hand,  causality
is lost at level $\eps >0$  and  is restored  just  in the limit $\eps \rightarrow 0$. Motivated by this fact, the
limit $\eps \rightarrow 0$ is usually referred to as \textit{causal
  limit}.

The WED approach contributes a variational approximation of the
limiting differential problem. As such, it paves the way to the
application of the general methods of the calculus of variations. The existence of the
approximations $u^\epsi$ is ascertained by the Direct
Method. This can be especially beneficial in actually computing
such approximations, for one can profit from the vast toolbox of
optimization methods. In addition, combinations with {\it other} approximations (additional parameters,
time and space discretization, data approximations) can be treated at
the $\eps>0$ level within the general frame of $\Gamma$-convergence~\cite{DalMaso93}.

The results of this paper can be compared with the ones of
\cite{Mielke2}. There, the dissipation is assumed to be independent of
the state and quadratic, namely, $\psi(v)=\|v\|_H^2/2$, corresponding
indeed to the classical gradient-flow case. The novelty here is in
considering a {\it state-dependent} dissipation instead, which has the
potential of considerably broaden the range of application of the
theory. This
extension has deep impact on the technical level, asking for the
development of a quite different analytical  approach.  Most notably,
the structure of the Euler--Lagrange equation features specific
state-dependent dissipative terms showing minimal integrability in
time. 
Note nonetheless that the simpler dissipation setting of \cite{Mielke2}
eventually allows more general
choices for the energy~$\phi$.

Before moving on, let us provide a brief review of the the vast
literature on the WED approach.  We however stress that all   
available results to date exclusively deal with the case $v\mapsto \psi(v)$,
 namely, with a dissipation which  is
independent of the state. 

In the gradient-flow  case of
$\psi$ quadratic, the use of the WED approach has to be traced back at
least to the study of Brakke mean-curvature flow of varifolds in
\cite{Ilmanen}.  Existence of periodic
solutions to gradient flows by the WED approach has also been obtained
\cite{Hirano94}. The general theory for the case of a $\lambda$-convex energy $\phi$ is in
\cite{Mielke2}. 
 Note that the linear case is also mentioned in the
classical PDE textbook by Evans \cite[Problem~3, p.~487]{Evans}. Nonconvex problems are discussed in  \cite{Akagi4},
see also \cite{Melchionna201sei} for nonpotential
perturbations,  \cite{LieroMelchionna} for
a stability result via $\Gamma$-convergence, and
  \cite{melchionna2} for an application to the study of
symmetries of solutions. Existence
of variational
solutions to the equation 
$
u_t -\nabla \cdot f(x,u,\nabla u) + \partial_u f(x,u,
\nabla u)=0
$
where the field $f$ is convex in
$(u, \nabla u)$ has been tackled in
\cite{Boegelein-et-al14}. Applications to stochastic PDE  
\cite{scarpastefanelli1,scarpastefanelli2} and to curves of maximal
slope in metric spaces
\cite{rossisavaresegattistefanelli1,rossisavaresegattistefanelli2} are
also available. We
also mention the applications to micro-structure
evolution 
 \cite{ContiOrtiz}, to  mean-curvature evolution
 of Cartesian surfaces \cite{spadarostefanelli}, and to
 the 
 incompressible Navier--Stokes system
 \cite{bathorystefanelli,ortizschmidtstefanelli}.

 The WED approach has been  also  considered  in the nonquadratic case,
 namely, for dissipations $v\mapsto \psi(v)$ of $q$ growth, for $1\leq
 q\not =2$. The rate-independent case $q=1$ \cite{mr} is tackled in
 \cite{MielkeOrtiz}, see also \cite{Mielke1} for a discrete
 counterpart. Applications to dynamic fracture \cite{deLucaDalMaso,LarsenOrtizRichardson} and dynamic plasticity
 are available \cite{DavoliStefanelli}.
 The doubly nonlinear setting $1<q\not =2$ has been
studied in \cite{AkagiMelchionnaStefanelli,Akagi1,Akagi2, Akagi3} under
different assumptions for the energy $\phi$, see also
\cite{AkagiMelchionna} for the case of nonpotential perturbations.

 The WED approach has been applied to hyperbolic and mixed
parabolic-hyperbolic problems, as well. Inspired by a  conjecture by De Giorgi \cite{DeGiorgi1996},
the case of semilinear waves $(\psi=0)$ is discussed in
\cite{serratilli1,stefanelli}. Extensions to other hyperbolic problems
\cite{Liero1,serratilli2}, also including forcings
\cite{tentarellitilli3,tentarellitilli,tentarellitilli2} or
dissipative terms $(\psi\not =0)$ \cite{marveggio,Liero2,serratilli2}, have already been
obtained.

The paper is organized as follows. In Section
\ref{notationpreliminaryresults}, we introduce notation, present some
preliminary  material and assumptions, and state our main result, namely, 
Theorem \ref{theorem}. Sections \ref{existenceeulerlagrange} and
\ref{causallimit} are devoted to the proof of the Theorem
\ref{theorem}. In particular, in Section \ref{existenceeulerlagrange}
we prove that the WED functional $W^{\eps}$ admits  global 
minimizers and that each  global  minimizer is a strong
solution of the Euler--Lagrange problem
\eqref{eulerolagrangeepsilon0}--\eqref{eulerolagrangeepsilon40}. 
Then,  in Section \ref{causallimit} we consider the causal limit
$\eps \rightarrow 0$ and we prove that the solutions of the
Euler--Lagrange problem
\eqref{eulerolagrangeepsilon0}--\eqref{eulerolagrangeepsilon40}
converge to the solutions of the gradient flow \eqref{eq:00}. 
Eventually, in Section \ref{sec:application} we present a PDE
example illustrating the abstract theory.

\section{ Main result}
\label{notationpreliminaryresults}
In this section, we introduce the  general functional setting 
 (Subsection \ref{sec:notation}),  record some preliminary
material (Subsection \ref{sec:preliminaries}), present our assumptions
(Subsection \ref{sec:assumptions}), and state our main result
(Subsection \ref{mainresult}). 

\subsection{Notation}\label{sec:notation}

 We assume  $X,\, H$  to  be two real and separable Hilbert spaces with $X\subset H$ densely and
compactly, so that, by indicating by $X^*$ the dual of $X$,
$(X,H,X^*)$ is a classical Hilbert triplet. We indicate by $( \cdot
, \cdot )$ the scalar product in $H$. The symbols $\lvert \cdot
\rvert_E$ and $\langle \cdot , \cdot \rangle_E$ will be used to
indicate the norm in the general Banach space $E$ and the duality
pairing between the dual $E^*$ and $E$, respectively.  

Given any convex, proper, and lower semicontinuous functional $f:H \to (-\infty,\infty]$, we
indicate by $D(f):=\{u\in H \ : \
f(u)<\infty\}$  its {\it effective domain} and  by  $\partial f (u)$ its {\it subdifferential} at a
point $u \in D(f)$ as the (possibly empty) set  
$$\partial f(u):= \{w \in H \ : \ (w,v-u) \leq f(v)- f(u) \  \ \forall v
\in H\}.$$
The map $\partial f: H \to 2^H$  (parts of $H$)   is a maximal monotone operator with
domain $D(\partial f ):= \{u\in D(f) \ : \ \partial f(u) \not =
\emptyset\}$. 
The {\it Legendre--Fenchel conjugate} $f^*$ of $f$ is  classically
  defined by
$f^*(w) := \sup_{u \in H} ( (w , u )- f(u))$ for $w \in H$ and is convex,
proper, and lower semicontinuous. Moreover, one has that $f^*(w) +
f(u) \geq (w, u )$ for all $u,\, w\in H$, as well as 
\begin{equation}
w \in \partial  f(u)  \ \Leftrightarrow \ u \in \partial  f^*(w) \
\Leftrightarrow \ f^*(w) + f(u) = (w , u).\label{moreau}
\end{equation}

For all convex, proper, and lower semicontinuous functionals $g:X \to (-\infty,\infty]$, we
define the subdifferential $\partial_X g: X\to 2^{X^*}$ at a point  
$u \in D(g)$ as the (possibly empty) set  
$$\partial_{X} g(u):= \{w \in X^* \ : \ \langle w,v-u\rangle_X \leq g(v)- g(u) \ \forall  v
\in X\}.$$

 Given a sequence $(f_\lambda)_\lambda$ and $f$ with $\, f_\lambda,\, f: H
\to [0,\infty)$ convex,  we say that
$f_\lambda$ {\it converge to $f$ in the sense of Mosco} \cite{attouch}
if
\begin{align*}
  \forall u_\lambda \rightharpoonup u \ \ \text{in} \ H \ \
  \Rightarrow \ \ f(u) \leq \liminf_{\lambda \to
  0}f_\lambda(u_\lambda), \\
  \forall v \in H \ \exists v_\lambda \to v: \quad \limsup_{\lambda
  \to 0} f_\lambda(v_\lambda) =f(v).
\end{align*}
In this context, $(v_\lambda)_\lambda$ is said to be a {\it recovery
  sequence} for $v$.

We indicate by $\d f: H \to H$ the  Fr\'echet differential of a
Fr\'echet-differentiable functional $f: H \to {\mathbb{R}}$. In the case of $\psi: H
\times H \to {\mathbb{R}}$, the symbols $\d_1$ and $\d_2$ denote the
partial   Fr\'echet differentials in the first and in the
second variable, respectively.

 \subsection{Preliminaries} \label{sec:preliminaries}

 We record here two preliminary lemmas, for later use.

\begin{lemma}[Aubin--Lions--Simon, {\cite[Thm.~3]{simon}}]\label{lemma:simon}
   The embedding
$$L^1(0,T;X) \cap  H^1  (0,T;H) \embed  C  ([0,T];H)$$ is compact.
\end{lemma}

\begin{lemma}[Integration by parts]
	\label{lemmaUlisse}
        Let $u\in H^1(0,T;H) \cap L^2(0,T;X)$ and $\xi \in L^2(0,T;H)$
        be such that $\dot{\xi} \in (H^1(0,T;H))^* +
        L^2(0,T;X^*)$. Let $t_1,t_2$ be Lebesgue points of $t \mapsto
        (\xi(t),u(t))$. Then, it holds that
	\begin{equation*}
		\langle  \dot{\xi} , u  \rangle_{H^1(t_1,t_2;H) \cap L^2(t_1,t_2;X)} = (\xi(t_2),u(t_2)) - (\xi(t_1),u(t_1)) - \int_{t_1}^{t_2} (\xi, \dot{u}).
	\end{equation*}
      \end{lemma}

       The latter   is an extension of \cite[Prop. 2.3]{Akagi2}.
The proof follows by adapting the argument  there, which was  originally developed  in
$L^r(t_1,t_2;X^*)$, and using the density of $L^r(t_1,t_2;X^*)$ into $
(H^1(t_1,t_2;X))^*$.

\subsection{ Assumptions}
\label{sec:assumptions}

 In addition to the general functional setting specified in
Subsection \ref{sec:notation}, we introduce here our assumptions on
data. We ask for the following.

\begin{itemize} 
	\item[\textbf{(A1)}] $\phi = \phi^1 + \phi^2 : H
          \rightarrow \mathbb{R}_+$, with $\phi^1(u)=\langle Au,u\rangle_X /2$
          if $u\in X$ and $\phi^1(u)=\infty$ if $u \in H\setminus X$
          for some $A:X \to X^*$ linear,
          continuous, and self-adjoint. There exist  $c_1, \,
          c_{2}>0 $  such that
		\begin{alignat}{4}
		\lvert u \rvert_X^2 &\leq c_{ 1} (1+\phi^1(u)  ) \;\; &&\forall u\in X ,\label{hypphy}\\
		\lvert Au \rvert_{X^*} &\leq c_{ 2} \lvert u \rvert_X \;\; &&\forall u\in X \label{ass4}.
		\end{alignat}
	The functional $\phi^2 : H \rightarrow \mathbb{R}_+$ is a
        proper, convex, and lower semicontinuous functional and there
        exists $c_{ 3}>0$ such that
	\begin{equation}
	\label{growthEta2Temp}
	\lvert \eta^2 \rvert_H^r \leq c_{ 3} (1+\phi^2(u)  ) \;\; \text{for every} \; u \in D(\partial \phi^2) \; \text{and} \; \eta^2 \in \partial \phi^2(u),
	\end{equation}
	for some $r \geq 2$.
	Moreover, $u_0 \in D(\phi)$.
\end{itemize}

\begin{itemize}
	\item[\textbf{(A2)}] $\psi \in C^2( H \times H ;
	\mathbb{R}_+)$, $\mathbb{R_+} := [0,\infty)$, is such that
	$\psi(u,\cdot)$ is convex and $\psi(u,0)=0$ for all $u \in
	H$. Moreover, there exists a positive constant $c_4$ such that
	\begin{equation}
	\lvert v \rvert_H^2 \leq c_4  \psi (u,v)   \quad \forall u,v \in H. \label{3ipotesinumero}
	\end{equation}	
	Furthermore, for all $R>0$, there exist positive constants
        $c_{5,R}$, $c_{6,R}$, $c_{7,R}$, $c_{8,R}$, $c_{9,R}$, and $c_{10,R}$ such that the following conditions hold:
	\begin{alignat}{3}
	&\lvert {\d}_1 \psi(u , v) \rvert_H \leq c_{5,R}
        \bigl(1+\lvert v \rvert_H^2  \bigr) \quad \forall u,v \in H \mbox{ such that } \lVert u \rVert_H \leq R, \label{0ipotesinumero}\\
	&\lvert {\d}_1 \psi(u_1 , v) - {\d}_1\psi(u_2,v) \rvert_H \leq
	c_{6,R} \lvert u_1 -u_2 \rvert_H \lvert v \rvert_H^2
	\nonumber\\
	&\qquad \qquad
	\forall u_1,u_2,v \in H \mbox{ such that } \lVert u_1 \rVert_H ,\, \lVert u_2 \rVert_H \leq R, \label{1ipotesinumero}\\
	&\lvert {\d}_2 \psi(u,v) \rvert_{H}^{2} \leq
        c_{7,R}\bigl(1+\lvert v \rvert_H^2  \bigr) \quad \forall u,v \in H \mbox{ such that }  \lVert u \rVert_H \leq R,\label{2ipotesinumero}\\
	&\|{\d}_2\psi(u_1,v)-{\d}_2\psi(u_2,v)\|_H \leq c_{8,R} \|
	u_1 - u_2 \|_H\|v\|_H \nonumber\\
	&\qquad \qquad \forall u_1, u_2, v \in H \mbox{ such that } \|u_1\|_H, \,  \|u_2\|_H\leq R,\label{23ipotesinumero}\\
	&\|{\rm d}_{21}\psi(u,v)w\|_H \leq c_{9,R}\|v\|_H \|w\|_H\quad \forall u, v, w \in H, \mbox{ such that } \|u\|_H\leq R,\label{21ipotesinumero}\\
	& ({\rm d}_{2 2}\psi(u,v)w,w) \geq c_{10,R} \| w \|_H^2  \quad \forall
	u, v,  w \in H, \mbox{ such that } \|u\|_H\leq
        R. \label{22ipotesinumero}
      \end{alignat}
      Here, we have used the notation $\d_{21}=\d_2\d_1$ and
      $\d_{22}=\d_2\d_2$.
    \end{itemize}
    

Before moving on, note that \eqref{3ipotesinumero} and
\eqref{2ipotesinumero} imply that  
\begin{align}
&\frac{1}{c_4}\|v\|_H^2 \leq  \psi(u,v) \leq ( {\d}_2\psi(u,v) , v )
                \leq c_{11,R} \bigl(1+\lvert v \rvert_H^2 ) \nonumber\\ &\qquad \qquad \forall u,v\in H \mbox{ such that }  \lVert u \rVert_H \leq R,\label{2.7ipotesinumero}
\end{align}
for $c_{11,R} = (c_{7,R} +1)/2$.  

By virtue of \eqref{23ipotesinumero} and \eqref{22ipotesinumero}, for
each $R > 0$, we can find positive constants $c_{12,R},
c_{13,R} > 0$ such that
\begin{align}\label{str-mono}
 \left( \d_2 \psi(u_1,v_1) - \d_2 \psi(u_2,v_2), v_1 - v_2 \right)
 \geq c_{12,R} \|v_1-v_2\|_H^2 - c_{13,R} \|u_1 - u_2\|_H^2 \|v_2\|_H^2
\end{align}
for $u_1,u_2,v_1,v_2 \in H$ satisfying $\|u_1\|, \|u_2\| \leq R$. Indeed, setting $v_\theta := (1-\theta)v_2 + \theta v_1$ for $\theta \in [0,1]$, we observe that
\begin{align*}
\MoveEqLeft{
\left( \d_2 \psi(u_1,v_1) - \d_2 \psi(u_2,v_2), v_1 - v_2 \right)
}\\[2mm]
&\geq \left( \d_2 \psi(u_1,v_1) - \d_2 \psi(u_1,v_2), v_1 - v_2 \right)
 + \left( \d_2 \psi(u_1,v_2) - \d_2 \psi(u_2,v_2), v_1 - v_2 \right)\\
&\geq \int^1_0 \left( \dfrac{\d}{\d \theta} \d_2 \psi(u_1,v_\theta), v_1 - v_2 \right) \, \d \theta - \| \d_2 \psi(u_1,v_2) - \d_2 \psi(u_2,v_2) \|_H \| v_1 - v_2 \|_H\\
&= \int^1_0 \left( \d_{22} \psi(u_1,v_\theta) (v_1 - v_2), v_1 - v_2 \right) \, \d \theta - \| \d_2 \psi(u_1,v_2) - \d_2 \psi(u_2,v_2) \|_H \|v_1-v_2\|_H\\
&\geq c_{10,R} \|v_1-v_2\|_H^2 - c_{8,R} \|u_1-u_2\|_H \|v_2\|_H \|v_1-v_2\|_H\\
&\geq \frac{c_{10,R}}2 \|v_1-v_2\|_H^2 - \frac{c_{8,R}^2}{2       c_{10,R}}  \|u_1-u_2\|_H^2 \|v_2\|_H^2
\end{align*}
so that \eqref{str-mono} follows by choosing $c_{12,R}=c_{10,R}/2$
and $c_{13,R}={c_{8,R}^2}/{(2c_{10,R})}$.

 We now prove a Mosco convergence result, which in particular
applies to the setting of Assumption (A2).

\begin{lemma}[Mosco convergence]\label{lem:2.3}  Let $\psi \in C^1(H \times H;\Rz_+)$
  fulfill \eqref{0ipotesinumero}. Then,  for all 
$u_\lambda \to u$ in $H$ one has that 
$
\psi(u_\lambda,  \cdot) \to \psi(u, \cdot)$  in the sense of Mosco  in $H$.
\end{lemma}

\begin{proof}
   Fix $u_1,u_2,v \in H$ and let $u_\theta := (1-\theta)u_1 +
  \theta u_2$ for $\theta \in (0,1)$. Owing to the fact that $\psi \in
  C^1(H \times H, \Rz_+)$, by using
  \eqref{0ipotesinumero} one computes 
  \begin{align*}
    &
\left| \psi(u_1,v) - \psi(u_2,v) \right|  = \int_0^1 \frac{\d}{\d
      \theta} \psi(u_\theta,v) \, \d \theta = \int_0^1
      (\d_1\psi(u_\theta,v),u_1 - u_2) \, \d \theta\\
    &\quad 
\leq \|u_1-u_2\|_H \int^1_0 \|\d_1 \psi(u_\theta,v)\|_H  \, \d \theta 
\leq c_{5,R}(1+\|v\|_H^2  ) \|u_1-u_2\|_H.
\end{align*}  This entails the pointwise convergence 
$
\psi(u_\lambda,v) \to \psi(u,v),
$
 so that, for all $v$ one can choose the constant recovery
sequence $v_\lambda=v$. 

 On the other hand,  let $v_\lambda  \rightharpoonup v$ in $H$. Then it follows that
\begin{align*}
 \psi(u_\lambda,v_\lambda) &\geq \psi(u,v_\lambda) - \left| \psi( u  ,v_\lambda) -
                 \psi( u_\lambda,v_\lambda)\right|\\
&\geq \psi(u,v_\lambda) - c_{5,R} (1+\|v_\lambda\|_H^2  ) \|u-u_\lambda\|_H,
\end{align*}
which implies that $
\liminf_{\lambda \to 0}  \psi(u_\lambda,v_\lambda)  \geq \psi(u,v)$.
\end{proof}

The argument of Lemma \ref{lem:2.3} can be easily generalized to
integral functionals on 
$L^2(0,T;H)$, also in presence of an arbitrary weight. More precisely, for all
nonnegative weights $\mu \in L^\infty(0,T)$, one also has that
\begin{equation}
\label{7ipotesinumero}
\int_0^T \mu \psi(u_k,\cdot) \rightarrow \int_0^T \mu  \psi(u,\cdot) \quad \text{in the sense of Mosco on }\,L^2(0,T;H)
\end{equation}
whenever $u_k \rightarrow u$ in $C([0,T];H)$. Note that in the
following we resort to the latter convergence by choosing
$\mu(t)=t$. 

Given Assumption (A1), we readily check that $\partial_X \phi^1|_X
\equiv A$ and $\partial \phi^1(u) = Au\in H$ for all $u \in D(\partial
\phi)=\{v\in X  \: : \: Au \in H\}$. Moreover, on account of
 assumptions \eqref{hypphy}--\eqref{ass4} and \eqref{2.7ipotesinumero}
 we can conclude that  $K(u_0)$ from \eqref{K} actually
corresponds to 
$$K(u_0) = \{u \in H^1(0,T;H)\cap L^2(0,T;X) \ : \ u(0)=u_0 \
\text{and} \ \phi^2(u)\in L^1(0,T)\}.$$

\subsection{Main result}
\label{mainresult}

Before stating our main result, i.e., Theorem \ref{theorem}, let us specify the 
Cauchy problem for the gradient flow \eqref{eq:00} as
\begin{align}
&\xi  + \eta^1 +\eta^2   = 0 \quad \text{in}\  H ,   \ \text{a.e. in} \
                (0,T) ,  \label{doublyNonLinearCauchy1}\\
&\xi  = {\d}_{2} \psi (u ,\dot{u} )\quad \text{in}\   H ,  \ \text{a.e. in} \
                (0,T) ,\label{doublyNonLinearCauchy2} \\
&\eta^1   = Au \quad\text{in} \  H  ,  \ \text{a.e. in} \
                (0,T), \label{doublyNonLinearCauchy3}\\
&\eta^2  \in \partial \phi^2(u )\quad\text{in} \  H ,  \ \text{a.e. in} \
                (0,T), \label{doublyNonLinearCauchybis}\\
&u(0)  = u_0.  \label{doublyNonLinearCauchy4}
\end{align}  
The {\it Euler--Lagrange} problem corresponding to the minimization of
$W^\eps$ given by \eqref{WEDfunctional} on $K(u_0)$ reads
\begin{align}
&-\eps\dot{\xi}_{\eps}  + \xi_{\eps}(t) +\eps\gamma_{\eps}  +
  \eta^1_{\eps}  + \eta^2_{\eps} = 0\ \text{in} \  X^* 
                                             ,   \ \text{a.e. in} \
                (0,T) ,\label{eulerolagrangeepsilon}
  \\
&\xi_{\eps}  =
       {\d}_{2}\psi(u_{\eps} ,\dot{u}_{\eps} )\quad\text{in}\  H  ,  \ \text{a.e. in} \
                (0,T) ,\\
&\gamma_{\eps}   =   {\d}_{1}\psi(u_{\eps}  ,\dot{u}_{\eps}       )\quad\text{in}\  H      ,  \ \text{a.e. in} \
                (0,T), \\
&\eta^1_{\eps}  = A u_{\eps} \ \text{in} \  X^* ,  \ \text{a.e. in} \
                (0,T) ,\\
&\eta^2_{\eps} \in \partial \phi^2(u_{\eps} )\ \text{in} \  H  ,  \ \text{a.e. in} \
                (0,T), \\
&u_{\eps}(0)=u_0 ,\\
& \epsi  \xi_{\eps}(T) = 0.\label{eulerolagrangeepsilon4}
\end{align} 
Note the occurrence of the natural Neumann condition $ \epsi  \xi_{\eps}(T) =
0$, which turns out to be instrumental for proving   {\it a priori}
estimates  and degenerates for $\epsi \to 0$. 
 
Our main result reads as follows.
\begin{theorem}[WED Approach]
	\label{theorem}
	Under Assumptions  {\rm (A1)--(A2)}   we have the
        following\/{\rm :} 
	\begin{enumerate}[{\rm (i)}]
		\item The {\rm WED} functional $W^{\eps}$ admits 
                  a   global minimizer  $u_{\eps} \in K(u_0)$\/{\rm ;}
		\item  Given  a    global  minimizer $u_{\eps}$ of $W^{\eps}$, the quintuple
                  $(u_{\eps},\xi_{\eps},\gamma_{\eps},\eta_{\eps}^1,\eta_{\eps}^2)$
                  with $\xi_{\eps} =
                  {\d}_2\psi(u_{\eps},\dot{u}_{\eps})$, $\gamma_{\eps}
                  = {\d}_1\psi(u_{\eps},\dot{u}_{\eps})$,
                  $\eta_{\eps}^1 = Au_{\eps}$,  and  $\eta_{\eps}^2 = \eps
                  \dot\xi_\eps - \xi_\eps-\eps\gamma_\eps-\eta^1_\eps\in \partial\phi^2(u_{\eps})$ belongs to
		\begin{align*}
		&[H^1(0,T;H)\cap L^2(0,T;X)] \times L^{2}(0,T;H) \times L^1(0,T;H) \\
		&\times L^2(0,T;X^*) \times L^{r}(0,T;H),
		\end{align*}
		and is a strong solution of the Euler--Lagrange problem \eqref{eulerolagrangeepsilon}--\eqref{eulerolagrangeepsilon4}\/{\rm ;}
		\item For any sequence $\eps_k \rightarrow 0$, there
                  exists a {\rm (}not relabeled{\rm )} subsequence such that \linebreak $(u_{\eps_k},\xi_{\eps_k},\eta^1_{\eps_k},\eta^2_{\eps_k}) \rightarrow (u,\xi,\eta^1,\eta^2)$ weakly* in
		\begin{equation*}
		[H^1(0,T;H) \cap L^2(0,T;X)] \times L^{2}(0,T;H) \times L^2(0,T;X^*) \times L^{r}(0,T;H),
		\end{equation*}
		$(u,\xi,\eta^1,\eta^2)$ is a strong solution of the
                gradient flow
                \eqref{doublyNonLinearCauchy1}--\eqref{doublyNonLinearCauchy4}
                 with  $\eta_1 = A u \in L^2(0,T;H)$.
	\end{enumerate}
\end{theorem}

The remainder of the paper is devoted to the proof of Theorem
\ref{theorem}. More precisely, parts ${\rm (i)}$ and ${\rm (ii)}$ are proved in
Section \ref{existenceeulerlagrange} by means of a regularization
procedure, and part ${\rm (iii)}$ is proved in Section \ref{causallimit}.

A {\it caveat} on notation: in the following we use the same symbol
$c$ to indicate a positive constant, possibly depending on $u_0, c_i, c_{i,R} ,
\ell$, but independent of $\eps$ and of the regularization parameter
$\lambda$ (introduced below). Note that the constant $c$ may change, even within the same line.

\section{Existence of solutions to the Euler--Lagrange problem}
\label{existenceeulerlagrange}
This section is devoted to the proof of Theorem \ref{theorem}.${\rm
  (i)}$-${\rm (ii)}$.  We show that $W^{\eps}$ admits global
minimizers $u_{\eps}$ and that such a global  minimizer  is a strong solution
of the Euler--Lagrange problem
\eqref{eulerolagrangeepsilon}--\eqref{eulerolagrangeepsilon4}. To this
aim, in Subsection \ref{directmethod} we define a regularized functional
$W^{\eps\lambda}$ for $\lambda > 0$ by replacing $\phi^1$ and
$\phi^2$ by their {\it Moreau--Yosida regularizations},  and check
that   $W^{\eps\lambda}$ admits  global  minimizers
$u_{\eps\lambda}$. These  $u_{\eps\lambda}$  are proved to
solve a  regularized  
Euler--Lagrange problem in Subsection \ref{regularizedeulerlagrange}. In
Subsection \ref{improvedintegrability}, we discuss  the integrability of
$u_{\eps\lambda}$. The chain rule for the potential $\psi$ is
addressed in Subsection \ref{chainrule}. This allows to derive in Subsection \ref{apriorilambda} some {\it a priori} estimates on the approximation
$u_{\eps\lambda}$ independently of $\lambda$. In Subsection
\ref{passagelimit}, we take the limit $\lambda \rightarrow 0$ (up to
subsequences) and prove that $u_{\eps} = \lim_{\lambda \rightarrow 0}
u_{\eps\lambda}$ solves the Euler--Lagrange problem
\eqref{eulerolagrangeepsilon}--\eqref{eulerolagrangeepsilon4}. Then, in
Subsection \ref{minimization}, we check that such $u_{\eps}$ minimizes
$W^{\eps}$, which proves Theorem  \ref{theorem}.${\rm (i)}$. Eventually,
we show that all  global  minimizers of $W^{\eps\lambda}$ solve the
Euler--Lagrange problem
\eqref{eulerolagrangeepsilon}--\eqref{eulerolagrangeepsilon4}, so that
Theorem \ref{theorem}.${\rm (ii)}$ holds, as well.

\subsection{The regularized functional}
\label{directmethod}
For every $\lambda >0$, let us introduce the regularized WED functional $W^{\eps\lambda}: L^2(0,T;H) \rightarrow [0,\infty]$ as
\begin{equation*}
W^{\eps\lambda}(u) =
\begin{cases}
\displaystyle\int_0^{T} e^{-t/\eps} \Big(
\eps\psi(J_{\lambda}u(t),\dot{u}(t)) + \phi^1_{\lambda}(u(t)) +
\phi_{\lambda}^2(u(t))  \Big) \;\d t  & \text{if}\; u\in K_0(u_0),\\[5mm]
\infty  & \text{else},
\end{cases}
\end{equation*}
where
$$
K_0(u_0) = \{u \in W^{1,2}(0,T;H) \ : \ u(0) = u_0\}.
$$
Here, $\phi^1_{\lambda} \in C^{1,1}(H;\mathbb{R})$ is the {\it
  Moreau--Yosida regularization} of $\phi^1$  in $H$,  namely
\begin{equation}
\label{defyosida}
\phi^1_{\lambda}(u) = \inf_{v\in H} \biggl( \frac{1}{2\lambda} \lvert u - v \rvert_H^2 + \phi^1(v)   \biggr) = \frac{1}{2\lambda} \lvert u - J_{\lambda}u \rvert_H^2 + \phi^1(J_{\lambda}u).
\end{equation}
The  map  $J_{\lambda} : H \to H$ is  the {\it resolvent}
of $\partial\phi^1$ at level $\lambda$, namely the linear solution
operator $J_{\lambda}: u  \in H  \mapsto J_{\lambda}u$ of the problem
\begin{equation*}
J_{\lambda}u - u + \lambda\partial\phi^1(J_{\lambda}u) = 0 \;\; \text{for every}\; u\in H.
\end{equation*}
We will also use the standard notation $A_\lambda:H \to H$ for the Yosida
approximation $A_\lambda u = (u - J_\lambda u)/\lambda = \d
\phi^1_\lambda (u)$ of $A$.
Analogously, we let $\phi_{\lambda}^2 \in C^{1,1}(H;\mathbb{R})$ be
the {Moreau--Yosida regularization} of $\phi^2$  in $H$  and $I_{\lambda}$ be
the corresponding resolvent of $\partial \phi^2$ at level $\lambda$.
\begin{lemma}[Global  Minimizers of the regularized functional]
  The regularized \emph{WED} functional
	$W^{\eps\lambda}$ admits  global  minimizers in $K_{0}(u_0)$.
	\begin{proof} 
The assertion follows by the Direct Method.  At first, note that
the constant-in-time function $u(t)\equiv u_0$ belongs to $K_0(u_0)$,
which is then not empty.  Let $( u_n)_n$ be an {\it
  infimizing sequence}, namely,  $u_n \in K_0(u)$ with   $W^{\eps\lambda}(u_n) \rightarrow
\inf_{ K_0(u_0)} W^{\eps\lambda}  $. As $W^{\eps\lambda}$ is coercive on $H^1(0,T;H)$ by assumption \eqref{3ipotesinumero}, we get
\begin{equation*}
\lVert u_n \rVert_{H^1(0,T;H)} \leq c.
\end{equation*}
Hence, up to not a relabeled subsequence, one has
\begin{equation}
\label{uinW1pLm}
u_n \rightharpoonup u \;\; \text{in} \; H^1(0,T;H).
\end{equation}
As $\phi^1_{\lambda}$ and $\phi_{\lambda}^2$ are  weakly  lower semicontinuous on $H$ we get
\begin{equation*}
\int_{0}^T e^{-t/\eps} \phi^j_{\lambda}(u) \leq \liminf_{n \rightarrow
  \infty} \int_{0}^T e^{-t/\eps} \phi^j_{\lambda}(u_n)\quad \text{for}
\ j=1,2.
\end{equation*}
Note that $u_0=u_n(0) \rightharpoonup u(0)$ in $H$, so that $u(0) = u_0$. This in particular implies that $u \in K_{0}(u_0)$. 
As $J_{\lambda}$ is a contraction in $H$, we obtain
\begin{equation*}
\biggl\lVert \frac{{\rm d}}{{\rm d}t}(J_{\lambda} u_n) \biggr\rVert_{L^2(0,T;H)} \leq \lVert \dot{u}_n \rVert_{L^2(0,T;H)} \leq c
\end{equation*}
and from the coercivity \eqref{hypphy} of $\phi^1$ we get
\begin{equation*}
\lVert J_{\lambda}u_n \rVert_{L^2(0,T;X)} \leq c.
\end{equation*}
Then, by applying the  Aubin--Lions--Simon Lemma \ref{lemma:simon},  again up to a not relabeled
subsequence,  we obtain the strong convergence
\begin{equation}
\label{aubin}
J_{\lambda}u_n \rightarrow J_{\lambda}u \;\; \text{in} \;
C([0,T];H).
\end{equation}
Moreover, since convergences \eqref{uinW1pLm}--\eqref{aubin} hold, relying on \eqref{7ipotesinumero} we have
\begin{equation*}
\int_{0}^T e^{-t/\eps}  \psi(J_{\lambda}u,\dot{u}) \leq \liminf_{n \rightarrow \infty} \int_{0}^T e^{-t/\eps}  \psi(J_{\lambda}u_n,\dot{u}_n).
\end{equation*}
This in turn guarantees that
\begin{equation*}
	W^{\eps\lambda} (u) \leq \liminf_{n \rightarrow \infty} W^{\eps\lambda}(u_n) = \inf_{K_{0}(u_0)} W^{\eps\lambda},
\end{equation*}
showing that $u$ is a  global  minimizer of $W^{\eps\lambda}$.
	\end{proof}
\end{lemma}
\subsection{Regularized Euler--Lagrange problem}
\label{regularizedeulerlagrange}
We now compute the first variation of the functional $W^{\eps\lambda}$
at a minimum point $u_{\eps\lambda} \in K_{0}(u_0)$. Fix $v \in  V := \{
v\in H^1(0,T;H) \; : \;v(0)= 0 \}$, so that $u_{\eps\lambda} + sv \in
K_{0}(u_0)$ for all $s \in \mathbb{R}$, and compute the derivative of the
real differentiable function $ s \mapsto  W^{\eps\lambda}(u_{\eps\lambda} + sv)$
at $s=0$. As $u_{\eps\lambda}$ is a  global  minimizer of $W^{\eps\lambda}$, one gets
\begin{align*}
0 &= \frac{\d}{\d s} \biggl[ W^{\eps\lambda}(u_{\eps\lambda} + sv)  \biggr]_{s=0}\\
 &=\int_{0}^T e^{-t/\eps} \frac{\d}{\d s} \biggl[ \eps\psi(J_{\lambda}(u_{\eps\lambda}+sv),\dot{u}_{\eps\lambda}+s\dot{v})  + \phi^1_{\lambda}(u_{\eps\lambda}+sv) + \phi_{\lambda}^2(u_{\eps\lambda} + sv)  \biggr]_{s=0}.
\end{align*}
Since $J_{\lambda}$ is linear we readily obtain
\begin{equation}
\label{temporaneo00}
0 = \int_{0}^T e^{-t/\eps} \biggl[   \eps (\xi_{\eps\lambda} , \dot{v} ) +\eps ( \gamma_{\eps\lambda} , J_{\lambda}v ) + ( \eta^1_{\eps\lambda} , v )  + ( \eta^2_{\eps\lambda} , v ) \biggr],
\end{equation}
where we have used the notation $\xi_{\eps\lambda} = {\d}_2\psi(J_{\lambda}u_{\eps\lambda} , \dot{u}_{\eps\lambda})$, $\gamma_{\eps\lambda} = {\d}_1\psi(J_{\lambda}u_{\eps\lambda} , \dot{u}_{\eps\lambda})$, $\eta^1_{\eps\lambda} = A_{\lambda}u_{\eps\lambda}$, and $\eta^2_{\eps\lambda} = \d\phi_{\lambda}^2(u_{\eps\lambda})$.
Integrating by parts the first term on the right-hand side of equation \eqref{temporaneo00} and recalling that $v$ vanishes at $t=0$, we get
\begin{equation*}
0= \int_{0}^T e^{-t/\eps} \biggl[  -\eps ( \dot{\xi}_{\eps\lambda} , v ) + ( \xi_{\eps\lambda} , v ) +\eps ( J_{\lambda}^*\gamma_{\eps\lambda} , v ) + ( \eta^1_{\eps\lambda} , v ) +  ( \eta^2_{\eps\lambda} , v )  \biggr] +\eps e^{-T/\eps}(\xi_{\eps\lambda}(T),v(T)) .
\end{equation*}
Here, we have again used the fact that $J_{\lambda}$ is linear and
continuous and we have indicated by $J_{\lambda}^*$ its adjoint. Since
the last integral equation holds for every $v \in  V$, by
further imposing $v(T) = 0$ we conclude that $u_{\eps\lambda}$ solves the regularized Euler--Lagrange problem
\begin{align}
&-\eps\dot{\xi}_{\eps\lambda}(t) + \xi_{\eps\lambda}(t) +\eps
  J_{\lambda}^*\gamma_{\eps\lambda}(t) + \eta^1_{\eps\lambda}(t) +
  \eta^2_{\eps\lambda}(t)=0\quad\text{in} \ 
                                 H,  \ \text{a.e. in} \ (0,T),\label{euleroepsilonlambda}\\ 
&\xi_{\eps\lambda}(t)  = {\d}_2\psi              (J_{\lambda}u_{\eps\lambda}(t),\dot{u}_{\eps\lambda}(t))\quad\text{in} \  H ,  \ \text{a.e. in} \ (0,T),\\
&\gamma_{\eps\lambda}(t)  = {\d}_1\psi (J_{\lambda}u_{\eps\lambda}(t),\dot{u}_{\eps\lambda}(t))\quad\text{in} \quad H ,  \ \text{a.e. in} \ (0,T),\\
&\eta^1_{\eps\lambda}(t)  = A_\lambda u_{\eps\lambda}(t)\quad\text{in}
                                                                                                                         \  H  \ \text{a.e. in} \ (0,T),\label{definizioneetalambda}\\
&\eta^2_{\eps\lambda}(t)  =       \d\phi_{\lambda}^2(u_{\eps\lambda}(t))\quad\text{in} \  H  \ \text{a.e. in} \ (0,T),\label{definizioneetalambda2}\\
&u_{\eps\lambda}(0)  = u_0 .
\end{align}
Eventually, by considering $v \in  V  $ with $v(T) \neq 0$, we
 additionally  deduce the terminal condition 
\begin{equation}
	 \eps \xi_{\eps\lambda}(T) = 0 \label{datoiniziale}.
      \end{equation}
      Note that the latter and 
      \eqref{2.7ipotesinumero} imply that $\dot u(T)=0$, as well.  

\subsection{Improved integrability}
\label{improvedintegrability}
We have proved that the regularized WED functional $W^{\eps\lambda}$
admits a  global  minimizer $u_{\eps\lambda} \in K_{0}(u_0)$ and that this solves
the regularized Euler--Lagrange problem
\eqref{euleroepsilonlambda}--\eqref{datoiniziale}. As $u_{\eps\lambda}$
belongs to $H^1(0,T;H)$, by recalling assumptions
\eqref{0ipotesinumero} and \eqref{2ipotesinumero} we have  that
 
$J_{\lambda}^* \gamma_{\eps\lambda} \in L^1(0,T;H)$,
$\xi_{\eps\lambda} \in L^{2}(0,T;H)$, and, by comparison in equation
\eqref{euleroepsilonlambda}, that $-\eps \dot{\xi}_{\eps\lambda} +
\eta^1_{\eps\lambda} + \eta_{\eps\lambda}^2 \in L^1(0,T;H)$. By using
the Lipschitz continuity of $ A_\lambda$   and $\d\phi^2_{\lambda}$ for $\lambda > 0$ fixed we can sharpen the statement on the integrability as follows.
\begin{lemma}[Improved integrability]\label{improved} Let
  $u_{\eps\lambda} \in K_0(u_0)$  minimize $W^{\eps\lambda}$. Then,
		\begin{align*}
	&u_{\eps\lambda} \in W^{ 2,\infty}(0,T;H),\\
	&\xi_{\eps\lambda} = {\d}_2 \psi(J_\lambda u_{\eps\lambda},\dot{u}_{\eps\lambda}) \in W^{1,\infty}(0,T;H),\\
	&\gamma_{\eps\lambda} = {\d}_1 \psi(J_\lambda u_{\eps\lambda},\dot{u}_{\eps\lambda}) \in L^{\infty}(0,T;H),\\
	& \eta^1_{\eps\lambda} =  A_\lambda  u_{\eps\lambda}  \in
   L^{\infty}(0,T;H), \\
   	& \eta^2_{\eps\lambda} =   \d \phi_{\lambda}^2 (u_{\eps\lambda}) \in L^{\infty}(0,T;H).
	\end{align*}
	\begin{proof}
As $A_\lambda ,\,  \d \phi_{\lambda}^2 : H \rightarrow H$ are Lipschitz continuous with Lipschitz constant $1/{\lambda}$, we have 
\begin{equation}
\label{boh}
	\lVert  \eta^1_{\epsi\lambda}  \rVert_{L^{\infty}(0,T;H)} \leq \frac{1}{\lambda} \lVert u_{\eps\lambda} \rVert_{L^{\infty}(0,T;H)} \leq c_{\lambda},
\end{equation}
and
\begin{equation*}
	\lVert  \eta^2_{\epsi\lambda} \rVert_{L^{\infty}(0,T;H)} \leq
        \frac{1}{\lambda} \lVert u_{\eps\lambda}
        \rVert_{L^{\infty}(0,T;H)} + c \leq c_{\lambda},
\end{equation*}
where, here and in the following, we explicitly record that the
constant $c_{\lambda}$ depends on $\lambda$. Given the fact that
$-\eps \dot{\xi}_{\eps\lambda} + \eta_{\eps\lambda}^1 +
\eta_{\eps\lambda}^2 = -\xi_{\eps\lambda}-\eps
J^*_\lambda\gamma_{\eps\lambda}  \in L^1(0,T;H)$, we deduce that $\dot{\xi}_{\eps\lambda} \in L^1(0,T;H)$. In particular,
	\begin{equation*}
	\lVert \xi_{\eps\lambda} \rVert_{W^{1,1}(0,T;H)} \leq c_{\lambda}.
	\end{equation*}
	 This,  along with \eqref{2.7ipotesinumero}, implies that $ u_{\eps\lambda} \in W^{1,\infty}(0,T;H)$ and ultimately
	\begin{equation}
	\label{stimaboh2}
	\lVert \gamma_{\eps\lambda} \rVert_{L^{\infty}(0,T;H)} \leq c_{\lambda}.
	\end{equation}
	 Then,  by using estimates \eqref{boh}--\eqref{stimaboh2}, and by a comparison in equation \eqref{euleroepsilonlambda} we get
	\begin{equation*}
	\lVert \xi_{\eps\lambda} \rVert_{W^{1,\infty}(0,T;H)} \leq
        \frac{c_{\lambda}}{ \eps}. 
      \end{equation*}
      
         As $\xi_{\eps\lambda}= {\rm d}_{2} \psi(J_\lambda u_{\epsi\lambda},\dot
		u_{\epsi\lambda})$, by using the
                smoothness of $\d_2\psi$ and differentiating in time one
                formally gets
		$$\dot\xi_{\eps\lambda}= {\rm d}_{21} \psi(J_\lambda u_{\epsi\lambda},\dot
		u_{\epsi\lambda}) J_\lambda \dot u_{\epsi\lambda}+ {\rm d}_{22} \psi(J_\lambda u_{\epsi\lambda},\dot
		u_{\epsi\lambda}) \ddot u_{\epsi\lambda}.$$
                The latter computation can indeed be made
                rigorous. In fact, one has that  $
		J_\lambda \dot u_{\epsi\lambda}\in L^\infty(0,T;H) $
                as $\dot u_{\epsi\lambda}\in L^\infty(0,T;H)  $ and
                $J_\lambda$ is Lipschitz continuous in $H$. Assumption
                \eqref{21ipotesinumero} entails that
                $$\|\d_{21}\psi(J_\lambda  u_{\epsi\lambda}, \dot
                u_{\epsi\lambda}) 	J_\lambda \dot
                u_{\epsi\lambda}\|_H\leq c_{9,R}\|\dot
                u_{\epsi\lambda} \|_H\| 	J_\lambda \dot
                u_{\epsi\lambda}\|_H,$$
                with $R= \| J_\lambda u_{\epsi\lambda} \|_{L^\infty(0,T;H)}$,
                hence $\d_{21}\psi(J_\lambda  u_{\epsi\lambda}, \dot
                u_{\epsi\lambda}) 	J_\lambda \dot
                u_{\epsi\lambda}\in L^\infty(0,T;H)$. By comparison
                $$   {\rm d}_{22} \psi(J_\lambda u_{\epsi\lambda},\dot
		u_{\epsi\lambda}) \ddot u_{\epsi\lambda}= \dot
                \xi_{\epsi\lambda}-\d_{21}\psi(J_\lambda  u_{\epsi\lambda}, \dot
                u_{\epsi\lambda}) 	J_\lambda \dot
                u_{\epsi\lambda}\in L^\infty(0,T;H)$$
                and one can use the coercivity of ${\rm d}_{22}
                \psi(J_\lambda u_{\epsi\lambda},\dot
                u_{\epsi\lambda})$, i.e., assumption
                \eqref{22ipotesinumero}, to get that $\ddot
                u_{\epsi\lambda}\in	L^\infty(0,T;H) $, as
                well. This proves the assertion.
                
		\end{proof}
\end{lemma}

\subsection{Chain rule for $\psi$}
\label{chainrule}

In view of deriving {\it a priori} estimates, it is paramount to prove a chain rule for $\psi$. Specifically, we have the following.
\begin{lemma}[Chain rule]
	\label{lemmachainrule}
	Let $u_{\eps\lambda}$ minimize $W^{\eps\lambda}$ in $K_{
          0}(u_0)$. Then, for all $t \in [0,T]$ the following holds
	\begin{align}
	\label{chainruleformula}
	\int_0^t \frac{\d}{\d s} \psi(J_{\lambda}u_{\eps\lambda} ,
          \dot{u}_{\eps\lambda}) =& \int_0^t \frac{\d}{\d s} (
                                    \xi_{\eps\lambda} ,
                                    \dot{u}_{\eps\lambda} ) + \int_0^t \biggl( -(  \dot{\xi}_{\eps\lambda} , \dot{u}_{\eps\lambda} ) + ( \gamma_{\eps\lambda} , J_{\lambda}\dot{u}_{\eps\lambda} ) \biggr).
	\end{align} 
	\begin{proof} As  $J_\lambda u_{\epsi\lambda},\, \dot
          u_{\epsi\lambda} \in W^{1,\infty}(0,T;H)$ from Lemma
          \ref{improved} and $\psi\in C^2(H\times H;\Rz_+)$, one has that $t
          \in (0,T)\mapsto \psi(J_\lambda u_{\epsi\lambda}(t), \dot
          u_{\epsi\lambda}(t))$ is differentiable and
          $$\frac{\d}{\d t}\psi(J_\lambda u_{\epsi\lambda}, \dot
          u_{\epsi\lambda}) = (\d_1 \psi(J_\lambda u_{\epsi\lambda}, \dot
          u_{\epsi\lambda}) ,J_\lambda \dot u_{\epsi\lambda})  + (\d_2 \psi(J_\lambda u_{\epsi\lambda}, \dot
          u_{\epsi\lambda}) ,\ddot u_{\epsi\lambda}).$$
More precisely, from
          \eqref{0ipotesinumero} and \eqref{2ipotesinumero} we deduce
          that
          \begin{align*}
            &\left|\frac{\d}{\d t}\psi(J_\lambda u_{\epsi\lambda}, \dot
          u_{\epsi\lambda}) \right| \leq \| \d_1 \psi(J_\lambda u_{\epsi\lambda}, \dot
          u_{\epsi\lambda})\|_H\| J_\lambda \dot u_{\epsi\lambda}\|_H  + \|\d_2 \psi(J_\lambda u_{\epsi\lambda}, \dot
              u_{\epsi\lambda}) \|_H\|\ddot u_{\epsi\lambda}\|_H\\
            &\quad \leq (c_{5,R}+c_{7,R})\big(1+\|\dot
              u_{\epsi\lambda}\|^2 \big) \big( \| J_\lambda \dot
              u_{\epsi\lambda}\|_H  +\|\ddot
              u_{\epsi\lambda}\|_H\big)
          \end{align*}
          with $R=\|J_\lambda u_{\epsi\lambda}\|_{L^\infty(0,T;H)}$,
          proving that  $t\mapsto \psi(J_\lambda u_{\epsi\lambda}(t), \dot
          u_{\epsi\lambda}(t)) \in W^{1,\infty}(0,T)$. In particular,
          the classical chain rule holds
		\begin{align*}
		\int_0^t \frac{\d}{\d s}
                  \psi(J_{\lambda}u_{\eps\lambda} ,
                  \dot{u}_{\eps\lambda}) = \int_0^t (
                  \xi_{\eps\lambda} , \ddot{u}_{\eps\lambda} ) +
                  \int_0^t ( \gamma_{\eps\lambda} ,
                  J_{\lambda}\dot{u}_{\eps\lambda} ) 
		\end{align*}
	and  equation \eqref{chainruleformula} readily follows.
	\end{proof}
\end{lemma}

\subsection{A priori estimates}
\label{apriorilambda}
In this section, we prove {\it a priori} estimates on the  global
 minimizers $u_{\eps\lambda}$  independently of  $\lambda$. To start
 with, test the regularized Euler--Lagrange equation
 \eqref{euleroepsilonlambda} by $\dot{u}_{\eps\lambda}$ and
 integrate it over $(0,T)$ to get
\begin{align*}
 \int_0^T \biggl( -\eps ( \dot{\xi}_{\eps\lambda} , \dot{u}_{\eps\lambda} ) + ( \xi_{\eps\lambda} , \dot{u}_{\eps\lambda} ) + \eps ( \gamma_{\eps\lambda}, J_{\lambda}\dot{u}_{\eps\lambda}) + ( \eta^1_{\eps\lambda} , \dot{u}_{\eps\lambda} ) + (\eta^2_{\eps\lambda} , \dot{u}_{\eps\lambda}) \biggr) = 0.
\end{align*}
 Owing to   the definitions \eqref{definizioneetalambda} and
\eqref{definizioneetalambda2} of $\eta^1_{\eps\lambda}$ and
$\eta^2_{\eps\lambda}$, using the chain rule
\eqref{chainruleformula} from Lemma \ref{lemmachainrule} with $t=T$
and recalling that $ \eps  \xi_{\eps\lambda}(T) =0$, we obtain 
\begin{align*}
   0=\;&\eps \int_0^T \frac{\d}{\d t} \psi(J_{\lambda}u_{\eps\lambda},\dot{u}_{\eps\lambda}) + \eps ( \xi_{\eps\lambda}(0) , \dot{u}_{\eps\lambda}(0) ) +  \int_0^T ( \xi_{\eps\lambda} , \dot{u}_{\eps\lambda} ) \\&+ \int_0^T \frac{\d}{\d t} \biggl(\phi^1_{\lambda}(u_{\eps\lambda}) + \phi_{\lambda}^2(u_{\eps\lambda})\biggr).
\end{align*}
By integrating the first and last terms and using $\dot{u}_{\eps\lambda}(T) = 0$ and $\psi(\cdot,0) = 0$, we infer that
\begin{align*}
0 =\; &\eps \biggl(  - \psi(J_{\lambda} u_0   ,\dot{u}_{\eps\lambda}(0)) +  ( \xi_{\eps\lambda}(0) , \dot{u}_{\eps\lambda}(0) ) \biggr) + \phi^1_{\lambda}(u_{\eps\lambda}(T)) - \phi^1_{\lambda}( u_0 ) \nonumber\\& + \phi_{\lambda}^2(u_{\eps\lambda}(T)) - \phi^2_{\lambda}( u_0 ) +  \int_0^T ( \xi_{\eps\lambda} , \dot{u}_{\eps\lambda} ).
\end{align*}
 This
ensures that  
\begin{equation}
\int_0^T ( \xi_{\eps\lambda} , \dot{u}_{\eps\lambda} ) +\phi^1_{\lambda}(u_{\eps\lambda}(T)) + \phi_{\lambda}^2(u_{\eps\lambda}(T))   \leq  \phi^1_{\lambda}( u_0 ) + \phi_{\lambda}^2( u_0 ).\label{5.23}
\end{equation}
Here, we also used the fact that $( \xi_{\eps\lambda}(0) ,
\dot{u}_{\eps\lambda}(0) ) -
\psi(J_{\lambda}u_0,\dot{u}_{\eps\lambda}(0)) =
\psi^*(J_{\lambda}u_0,\xi_{\eps\lambda}(0))  \geq 0$ due to the
identity  \eqref{moreau} applied to $v \mapsto \psi(J_\lambda
u_0,v)$, as well as $\psi(J_\lambda u_0,0) = 0$.  By means  inequality \eqref{2.7ipotesinumero} and  
the nonnegativity of $\phi^1_{\lambda}$ and $\phi_{\lambda}^2$ we deduce that
\begin{equation}
\label{154}
\int_0^T \lvert \dot{u}_{\eps\lambda} \rvert_H^2 \leq  c.
\end{equation}
Via \eqref{0ipotesinumero} and \eqref{2ipotesinumero} this implies that
\begin{equation}
\label{155}
\int_0^T \lvert \xi_{\eps\lambda}\rvert_{H}^{2} \leq c,
\end{equation}
and
\begin{equation}
\label{157}
\int_0^T \lvert \gamma_{\eps\lambda} \rvert_{H} \leq c.
\end{equation}

We test again \eqref{euleroepsilonlambda} by  $\dot{u}_{\eps\lambda}$
and integrate it on $(0,t)$ with $t\in(0,T)$. Again from the chain rule
\eqref{chainruleformula} we obtain
\begin{align*}
&\phi^1_{\lambda}(u_{\eps\lambda}(t)) + \phi_{\lambda}^2
                 (u_{\eps\lambda}(t)) + \int_0^t  ( \xi_{\eps\lambda}
                 , \dot{u}_{\eps\lambda} ) + \eps \Big( (
                 \xi_{\eps\lambda}(0) , \dot{u}_{\eps\lambda}(0) ) -
                 \psi(J_{\lambda} u_0   , \dot{u}_{\eps\lambda}(0)) 	\Big) =\\ 
&=\eps \Big( ( \xi_{\eps\lambda}(t) , \dot{u}_{\eps\lambda}(t) ) -\psi(J_{\lambda}u_{\eps\lambda}(t) , \dot{u}_{\eps\lambda}(t))	\Big) +\phi^1_{\lambda}( u_0  ) + \phi_{\lambda}^2( u_0  ).
\end{align*}
Thanks to \eqref{2.7ipotesinumero}  and the fact that $\eps ( (
                 \xi_{\eps\lambda}(0) , \dot{u}_{\eps\lambda}(0) ) -
                 \psi(J_{\lambda} u_0   ,
                 \dot{u}_{\eps\lambda}(0)) 	)\geq0$  we deduce that
\begin{align*}
  \phi^1_{\lambda}(u_{\eps\lambda}(t)) +
    \phi_{\lambda}^2(u_{\eps\lambda}(t)) &\leq  \eps \biggl( (
    \xi_{\eps\lambda}(t) , \dot{u}_{\eps\lambda}(t) )
    -\psi(J_{\lambda}u_{\eps\lambda}(t) , \dot{u}_{\eps\lambda}(t))
    \biggr)\\ &\quad +\phi^1_{\lambda}( u_0  ) + \phi_{\lambda}^2( u_0 ).
\end{align*}
Integrating it once more on $(0,T)$ we get
\begin{equation}
\label{202}
\int_0^T \phi^1_{\lambda}(u_{\eps\lambda}) + \int_0^T
\phi_{\lambda}^2(u_{\eps\lambda})  \leq c\eps + T
\phi^1_{\lambda}(u_0) +  T  \phi_{\lambda}^2(u_0)  \leq  c.
\end{equation}
From the bound on $\phi_{\lambda}^1(u_{\eps\lambda})$, by using assumption \eqref{hypphy}, we obtain
\begin{equation}
\label{200}
\int_0^T \lvert J_{\lambda}u_{\eps\lambda} \rvert_{X}^2 \leq c.
\end{equation}
Hence, thanks to assumption \eqref{ass4}, we also get
\begin{equation*}
\int_0^T \lvert A_\lambda u_{\eps\lambda} \rvert_{X^*}^2 \leq c.
\end{equation*}
From the bound on $\phi_{\lambda}^2(u_{\eps\lambda})$, by using $\eta^2_{\eps\lambda} = \partial
\phi^2_\lambda(u_{\varepsilon\lambda}) \in \partial \phi^2(I_\lambda
u_{\varepsilon\lambda})$, $\phi^2(I_\lambda  u_{\epsi\lambda} ) \leq \phi_\lambda^2( u_{\epsi\lambda} )$, and assumption \eqref{growthEta2Temp}, we have
\begin{equation}
	\int_0^T \lvert \eta_{\eps\lambda}^2 \rvert_H^r \leq c.
\end{equation}
Eventually, by comparison in the regularized Euler--Lagrange equation \eqref{euleroepsilonlambda}, we get
\begin{align}
\label{201}
\varepsilon\lVert \dot{\xi}_{\eps\lambda} \rVert_{L^1(0,T;H) + L^2(0,T;X^*)} \leq c.
\end{align}

\subsection{Passage to the limit as $\lambda \rightarrow 0$}
\label{passagelimit}
We aim at checking that   by taking the limit of
$(u_{\eps\lambda},\xi_{\eps\lambda},\gamma_{\eps\lambda},\eta^1_{\eps\lambda},\eta^2_{\eps\lambda})$
as $\lambda \to 0$ we obtain a   strong solution   of the Euler--Lagrange problem \eqref{eulerolagrangeepsilon}--\eqref{eulerolagrangeepsilon4}. Thanks to estimates \eqref{154}--\eqref{157} and \eqref{200}--\eqref{201} we find $(u_{\eps},\xi_{\eps},\gamma_{\eps},\eta^1_{\eps},\eta^2_{\eps})$ such that, up to not relabeled subsequences, the following convergences hold\/{\rm :}
\begin{alignat}{4}
u_{\eps\lambda} &\rightharpoonup u_{\eps}  \quad &&\text{in} \quad H^1(0,T;H),\label{convergencelambdau}\\
J_{\lambda}u_{\eps\lambda} &\rightharpoonup \beta_{\eps}
                                                         \quad &&\text{in}
                                                         \quad   H^1(0,T;H)     \cap L^2(0,T;X),\label{convergencelambdajux}\\
J_{\lambda}^*\gamma_{\eps\lambda} &  \rightharpoonup   \gamma_{\eps} \quad &&\text{in} \quad (H^1(0,T;H))^*\label{convergencelambdagamma},\\
\gamma_{\eps\lambda} &  \rightharpoonup \tilde{\gamma}_{\eps} \quad &&\text{in} \quad (H^1(0,T;H))^*\label{convergencelambdajgamma},\\
A_\lambda u_{\eps\lambda} &\rightharpoonup \eta_{\eps}^1 \quad &&\text{in} \quad L^2(0,T;X^*),\label{convergencelambdaeta}\\
\eta^2_{\eps\lambda} &\rightharpoonup \eta^2_{\eps} \quad &&\text{in} \quad L^{r}(0,T;H), \label{convergenceeta2lambdaTemp}\\
\xi_{\eps\lambda} &\weakstar \xi_{\eps} \quad&&\text{in} \quad L^2(0,T;H) \cap  BV([0,T];X^*),\label{convergencelambdaxi}\\
\dot{\xi}_{\eps\lambda} &\rightharpoonup \dot{\xi}_{\eps}  \quad &&\text{in} \quad (H^1(0,T;H))^* + L^2(0,T;X^*).\label{convergencelambdaxipunto}
\end{alignat}
Convergences \eqref{convergencelambdagamma},
\eqref{convergencelambdajgamma}, and \eqref{convergencelambdaxipunto}
are obtained from the estimate \eqref{157} by considering the
continuous embedding of $L^1(0,T;H)$ into the Hilbert space
$(H^1(0,T;H))^*$ and by noting that $\|J_\lambda^*\|_{\mathcal{L}(H)}
\leq 1$, where $\|\,\cdot\,\|_{\mathcal{L}(H)}$ denotes the 
operator  norm for bounded linear operators in $H$. Thanks to the
 Aubin--Lions--Simon Lemma \ref{lemma:simon},  from convergence
\eqref{convergencelambdajux} we deduce the strong convergence 
\begin{equation}
\label{convergencestrongjuv}
J_{\lambda}u_{\eps\lambda} \rightarrow \beta_{\eps} \;\; \text{in} \; C([0,T];H).
\end{equation}
 We  can pass to the limit into the regularized Euler--Lagrange equation \eqref{euleroepsilonlambda} and obtain
\begin{equation}
\label{limitequationepsilon}
-\eps \dot{\xi}_{\eps} + \xi_{\eps} +\eps\gamma_{\eps} + \eta^1_{\eps} + \eta^2_{\eps} = 0 \;\; \text{in} \;\; (H^1(0,T;H))^* + L^2(0,T;X^*).
\end{equation}

Let us now identify the limit functions. We proceed in subsequent
steps.
	 
$\bullet$	 {\it Step 1: Identification $\tilde{\gamma}_{\eps} = \gamma_{\eps}$.}
	By means of convergence \eqref{convergencelambdagamma} we have, for all $v\in H^1(0,T;H)$,
	\begin{equation*}
	 \langle J_{\lambda}^* \gamma_{\eps\lambda} , v  \rangle_{H^1(0,T;H)}\rightarrow
	  \langle \gamma_{\eps} , v  \rangle_{H^1(0,T;H)}.
	\end{equation*} 
	On the other hand, as $J_{\lambda}$ is linear, we get
	\begin{equation*}
	 \langle J_{\lambda}^* \gamma_{\eps\lambda} , v  \rangle_{H^1(0,T;H)} =  \langle \gamma_{\eps\lambda} , J_{\lambda}v  \rangle_{H^1(0,T;H)} \rightarrow \langle \tilde{\gamma}_{\eps} , v \rangle_{H^1(0,T;H)} ,
	\end{equation*}
	where we have used the strong convergence of $J_{\lambda}v$ to
        $v$ in $H^1(0,T;H)$ and the weak convergence
        \eqref{convergencelambdajgamma} of $\gamma_{\eps\lambda}$ to
        $\tilde{\gamma}_{\eps}$ in $(H^1(0,T;H))^*$. Hence, we have $\tilde{\gamma}_{\eps} = \gamma_{\eps} $.
	
$\bullet$	 {\it Step 2: Strong convergence of $u_{\eps\lambda}$ to $u_{\eps}$.}
	Recalling the definition of the Moreau--Yosida regularization \eqref{defyosida} and bound \eqref{202} on $\phi^1_{\lambda}(u_{\eps\lambda})$ in $L^1(0,T;\mathbb{R}_+)$, we get the following strong convergence\/{\rm :}
	\begin{equation}
	\label{convergencetemp}
	\lVert u_{\eps\lambda} - J_{\lambda}u_{\eps\lambda} \rVert_{L^2(0,T;H)}^{2} \leq 2\lambda \lVert \phi^1_{\lambda}(u_{\eps\lambda}) \rVert_{L^1(0,T)} \rightarrow 0 .
	\end{equation}
	Let $\varphi \in L^2(0,T;H)$ be given. Using the weak convergence \eqref{convergencelambdau} of $u_{\eps\lambda}$ in $H^1(0,T;H)$ along with the strong convergence \eqref{convergencetemp}, we infer that $J_{\lambda}u_{\eps\lambda}$ weakly converges to $u_{\eps}$ in $L^2(0,T;H)$, namely,
	\begin{align*}
	\biggl| \int_{0}^T  (& \varphi , J_{\lambda}u_{\eps\lambda} -
                               u_{\eps} ) \biggr| \leq \int_{0}^T | (
                               \varphi , J_{\lambda}u_{\eps\lambda} -
                               u_{\eps\lambda} ) | +\left| \int_0^T  ( \varphi , u_{\eps\lambda} - u_{\eps} ) \right| \to 0.
	\end{align*}
	Hence, thanks to the uniqueness of the limit and to strong convergence \eqref{convergencestrongjuv} of $J_{\lambda}u_{\eps\lambda}$ to $\beta_{\eps}$ in $C([0,T];H)$, we deduce that
	\begin{equation}
	\label{convergencejuc0}
	J_{\lambda}u_{\eps\lambda} \rightarrow u_{\eps} \;\; \text{in} \;  C([0,T];H).
	\end{equation}
	From the strong convergences \eqref{convergencetemp} and \eqref{convergencejuc0} we deduce that
	\begin{equation*}
	\lVert u_{\eps\lambda} - u_{\eps} \rVert_{L^2(0,T;H)} \leq \lVert u_{\eps\lambda} - J_{\lambda}u_{\eps\lambda} \rVert_{L^2(0,T;H)} +  \lVert J_{\lambda}u_{\eps\lambda} - u_{\eps} \rVert_{L^2(0,T;H)} \rightarrow 0,
	\end{equation*}
	which implies the pointwise convergence, up to a not relabeled subsequence,
	\begin{equation}
	u_{\eps\lambda}(t) \rightarrow u_{\eps}(t) \;\; \text{in} \; H, \; \text{for a.e.} \;t\in (0,T).\label{nome}
	\end{equation}
	This pointwise convergence and bound \eqref{154} on $u_{\eps\lambda}$ allow us to apply the Vitali--Lebesgue Theorem (or H\"older's inequality) and obtain the strong convergence
	\begin{equation}
	\label{convergenceulq}
	u_{\eps\lambda} \rightarrow u_{\eps} \;\; \text{in} \; L^q(0,T;H)\;\;\text{for any}\; q<\infty.
	\end{equation}

$\bullet$	 {\it Step 3: Identification of $\eta^1_{\eps}$}.
	As $A_\lambda u_{\eps\lambda} =
        A(J_{\lambda}u_{\eps\lambda})$, $A$ is linear and continuous from $X$ to
        $X^*$, and  
        $J_{\lambda}u_{\eps\lambda}$ weakly converges to $u_{\eps}$ in
        $L^2(0,T;X)$ from \eqref{convergencelambdajux} and
        \eqref{convergencejuc0}, it follows that
	\begin{equation}
	\label{etaespilonlambdaconverge}
	A(J_{\lambda}u_{\eps\lambda}) \rightharpoonup A u_{\eps} \;\; \text{in}\;L^2(0,T;X^*).
	\end{equation}
In particular,
	\begin{equation*}
	\eta^1_{\eps}  = A u_{\eps}  \;\; \text{in} \; X^*,
        \;\text{a.e. in} \  (0,T).
	\end{equation*}
	
$\bullet$	{\it Step 4: Identification of $\eta^2_{\eps}$.}
	 From   the definition \eqref{defyosida} of the Moreau--Yosida regularization and bound \eqref{202} on $\phi_{\lambda}^2(u_{\eps\lambda})$ in $L^1(0,T)$ we get
	\begin{equation*}
		\lVert u_{\eps\lambda} - I_{\lambda}u_{\eps\lambda} \rVert_{L^2(0,T;H)}^{2} \leq 2 \lambda \lVert \phi_{\lambda}^2(u_{\eps\lambda}) \rVert_{L^1(0,T)} \rightarrow 0,
	\end{equation*}
where  we recall that  $I_\lambda = (I + \lambda \partial
\phi^2)^{-1}$ denotes the resolvent for $\partial
\phi_\lambda^2$.  This,  along with the strong convergence \eqref{convergenceulq}, implies
	\begin{equation}
	\label{convergenzaTemp}
	I_{\lambda}u_{\eps\lambda} \rightarrow u_{\eps} \;\; \text{in} \; L^2(0,T;H).
	\end{equation}
	From the convexity of $\phi^2$, thanks to the strong
        convergence \eqref{convergenzaTemp} and to the weak
        convergence \eqref{convergenceeta2lambdaTemp}, implies that
	\begin{equation*}
	\label{convergenceeta2Temp}
	\int_0^T  (\eta^2_{\eps\lambda} , I_{\lambda} u_{\eps\lambda}) \rightarrow  \int_0^T (\eta^2_{\eps} , u_{\eps}),
      \end{equation*}
       so that the identification   
	\begin{equation*}
	\eta^2_{\eps} \in \partial\phi^2(u_{\eps} ) \;\; \text{in} \;
        H, \; \text{a.e. in} \ (0,T)
      \end{equation*}
       follows by   \cite[Prop.~2.5, p.~27]{Brezis73} .

$\bullet$	 {\it Step 5: Identification of $\xi_{\eps}$.}
                
As $\xi_{\eps\lambda} \in W^{1,\infty}(0,T;H)$
        and $u_{\eps\lambda} \in H^1(0,T;H)$, we have that $t\mapsto (
        \xi_{\eps\lambda}(t) , u_{\eps\lambda}(t) ) \in
        H^1(0,T)$. Also using \eqref{datoiniziale} we can
        integrate by parts and get 
	\begin{align*}
	\int_{t}^{T} ( \eps \xi_{\eps\lambda} , \dot{u}_{\eps\lambda}
          ) =\;&  - \eps ( \xi_{\eps\lambda}(t) , u_{\eps\lambda}(t) )
                 + \int_{t}^{T}  ( -\eps \dot{\xi}_{\eps\lambda} ,
                 u_{\eps\lambda} ) \quad \forall t \in [0,T].
	\end{align*}
        We now rewrite the term in  $\eps\dot{\xi}_{\eps\lambda}$
        by
using the regularized Euler--Lagrange equation
\eqref{euleroepsilonlambda} and integrate again over $(0,T)$
getting 
	\begin{align}
	\label{l1temp}
	&\int_0^T \!\!\int_{t}^{T} ( \eps \xi_{\eps\lambda} ,
          \dot{u}_{\eps\lambda} )= - \eps \int_0^T ( \xi_{\eps\lambda} ,
   u_{\eps\lambda} )  - \varepsilon \int_0^T  \!\!\int_{t}^{T} (
   J_{\lambda}^* \gamma_{\eps\lambda}, u_{\eps\lambda})\nonumber
          \\
          &\quad  - \int_0^T \!\!\int_{t}^{T} ( \xi_{\eps\lambda} , u_{\eps\lambda} )- \int_0^T  \!\!\int_{t}^{T} ( \eta^2_{\eps\lambda}, u_{\eps\lambda} ) - \int_0^T \!\!\int_{t}^{T} ( \eta^1_{\eps\lambda}, u_{\eps\lambda} ).
	\end{align}

        Our aim is to pass to the $\limsup$ as $\lambda \to 0$ in
         \eqref{l1temp}.	Let us consider each term in the right-hand side separately. By means of the strong convergence \eqref{convergenceulq} of $u_{\eps\lambda}$ to $u_{\eps}$ in $L^q(0,T;H)$ for any $q<\infty$ and of the weak convergence \eqref{convergencelambdaxi} of $\xi_{\eps\lambda}$ to $\xi_{\eps}$ in $L^{2}(0,T;H)$ we obtain
	\begin{equation}
	\label{convlebesguepoints}
	\int^T_{t}( \xi_{\eps\lambda} , u_{\eps\lambda} )
        \rightarrow \int^T_{t} ( \xi_{\eps} , u_{\eps} )
        \quad \forall t \in [0,T],
      \end{equation}
	as well as
	\begin{equation}
	\label{nonloso}
	 \int_0^T\!\! \int_{t}^{T} ( \xi_{\eps\lambda} , u_{\eps\lambda} ) \rightarrow \int_0^T\!\! \int_{t}^{T} ( \xi_{\eps} , u_{\eps} ). 
      \end{equation}
      
The $J^*_\lambda \gamma_{\eps\lambda}$-term on the
right-hand side of \eqref{l1temp} can be treated as follows
	\begin{align*}
&\int_0^T\!\!	\int_{t}^{T} (   J_{\lambda}^*\gamma_{\eps\lambda} , u_{\eps\lambda} ) =
          \int_0^T t \, (   J_{\lambda}^*\gamma_{\eps\lambda} ,
          u_{\eps\lambda} )   =  \int^T_{0}(
          \gamma_{\eps\lambda} , t J_{\lambda}u_{\eps\lambda} ) 
          \\
          &\quad = \int^T_{0}( \gamma_{\eps\lambda},  t J_{\lambda}
           u_{\eps\lambda} -  t  u_{\eps} ) +  \int^T_{0
           } ( \gamma_{\eps\lambda} ,  t u_{\eps} ).
	\end{align*}
   The  strong convergence \eqref{convergencejuc0} of
   $J_{\lambda}u_{\eps\lambda}$ to $u_{\eps}$ in $C([0,T];H)$ and
   bound \eqref{157} of $\gamma_{\eps\lambda}$ in $L^1(0,T;H)$
   ensure that  the first integral in the above right-hand side converges to $0$
as $\lambda \to 0$. Moreover, the weak convergence
\eqref{convergencelambdajgamma} of $\gamma_{\eps\lambda}$ to
$\gamma_{\eps}$ in $(H^1(0,T;H))^*$ ensures that 
$$\int_0^T (\gamma_{\eps\lambda},tu_\eps ) \to  \langle \gamma_{\eps}
         , tu_{\eps} \rangle _{H^1(0,T;H)}.$$
         All in all, we have proved that
         \begin{equation}
	\label{convergenzagammau}
\int_0^T\!\! 	\int_{t}^{T} (   J_{\lambda}^*\gamma_{\eps\lambda} ,
u_{\eps\lambda} ) \to  \langle \gamma_{\eps}
         , tu_{\eps} \rangle _{H^1(0,T;H)}.
       \end{equation}

	From convergences \eqref{convergenceulq} and \eqref{convergenceeta2lambdaTemp} we get
	\begin{equation*}
		\int_t^T (\eta_{\eps\lambda}^2 , u_{\eps\lambda})
                \rightarrow \int_t^T (\eta_{\eps}^2 , u_{\eps}) \quad
                 \forall t \in [0,T]
              \end{equation*}
      as well as
              	\begin{equation}
	\label{convTempo}
		\int_0^T\!\!\int_t^T (\eta_{\eps\lambda}^2 , u_{\eps\lambda})
                \rightarrow \int_0^T\!\!\int_t^T (\eta_{\eps}^2 ,
                u_{\eps}).
              \end{equation}
              
	As for the last term in  the right-hand side of 
        \eqref{l1temp}, for all $t\in [0,T]$ we use  the fact that  
	\begin{align*}
	\int_{t}^{T}  \langle A_\lambda u_{\eps\lambda},
          J_{\lambda}u_{\eps\lambda} \rangle_X &=  \int_{t}^{T} (
                                               A_\lambda
                                               u_{\eps\lambda},
                                               u_{\eps\lambda} )
                                               -\lambda
                                               \int_{t}^{T} \lvert
                                               A_\lambda
                                               u_{\eps\lambda} \rvert_{H}^2\leq \int_{t}^{T} ( A_\lambda u_{\eps\lambda}, u_{\eps\lambda} ),
	\end{align*}
	and rewrite it as
	\begin{equation}
	\label{ineqtemp}
		- \int_{t}^{T} ( A_\lambda u_{\eps\lambda}, u_{\eps\lambda} ) \leq - \int_{t}^{T}  \langle A_\lambda u_{\eps\lambda}, J_{\lambda}u_{\eps\lambda} \rangle_X.
	\end{equation}
	Since $A$ is maximal monotone and
        $A_\lambda u_{\eps\lambda} = A(J_\lambda u_{\eps\lambda})$, we have
	\begin{equation*}
	0 \leq \int_{t}^{T} \langle A_\lambda u_{\eps\lambda}- A u_\eps , J_{\lambda}u_{\eps\lambda} - u_{\eps} \rangle_X,
	\end{equation*} 
	from which, thanks to the weak convergences
        \eqref{etaespilonlambdaconverge} and
        \eqref{convergencelambdajux} of $A_\lambda u_{\eps\lambda}$ to
        $A u_\eps$ in \linebreak $L^2(0,T;X^*)$ and of $J_{\lambda}u_{\eps\lambda}$ to $u_{\eps}$ in $L^2(0,T;X)$, respectively, it follows that
	\begin{equation*}
	\int_0^T\!\! \int_{t}^{T}  \langle A u_\eps , u_{\eps} \rangle_X \leq \liminf_{\lambda\rightarrow 0} \int_{0}^{T}\!\!  \int_t^T  \langle A_\lambda u_{\eps\lambda}, J_{\lambda}u_{\eps\lambda} \rangle_X,
	\end{equation*}
	that is
	\begin{equation}
	\label{limsupphi}
	\limsup_{\lambda\rightarrow 0} \biggl( - \int_0^T\!\! \int_{t}^{T} \langle A_\lambda u_{\eps\lambda}, J_{\lambda}u_{\eps\lambda} \rangle_X  \biggr) \leq -\int_0^T\!\! \int_{t}^{T}  \langle A u_\eps, u_{\eps} \rangle_X.
	\end{equation}
	Eventually, by combining inequalities \eqref{ineqtemp} and
        \eqref{limsupphi}  we get
	\begin{equation}
	\label{limsupinequatemp}
		\limsup_{\lambda\rightarrow 0} \biggl(  - \int_0^T\!\!  \int_{t}^{T} ( A_\lambda u_{\eps\lambda}, u_{\eps\lambda} )    \biggr) \leq -\int_0^T\!\! \int_{t}^{T}  \langle A u_\eps, u_{\eps} \rangle_X.
	\end{equation}
	
	We now pass to the $\limsup$ as $\lambda \rightarrow 0$, in equation \eqref{l1temp}. Using convergences \eqref{convlebesguepoints}--\eqref{convTempo} and inequality \eqref{limsupinequatemp} we obtain
	\begin{align*}
	\limsup_{\lambda \rightarrow 0}&  \int_0^T\!\! \int_{t}^{T} (
                                         \eps \xi_{\eps\lambda} ,
                                         \dot{u}_{\eps\lambda} ) \leq
                                         \;  - \eps \int_0^T (
                                         \xi_{\eps} , u_{\eps} ) -
                                          \langle \gamma_{\eps}, t
                                         u_{\eps} \rangle_{H^1( 0,T;H)}  \\ & -  \int_0^T\!\! \int_{t}^{T}  \langle A u_\eps, u_{\eps} \rangle_X - \int_0^T\!\! \int_{t}^{T} (\eta^2_{\eps} , u_{\eps}) - \int_0^T\!\!  \int_{t}^{T} ( \xi_{\eps} , u_{\eps} ).
	\end{align*}
	By using the limit equation \eqref{limitequationepsilon} 
        and the fact that $\eta^1_\eps = Au_\eps$, we rewrite the
        above right-hand side as  
	\begin{align*}
		\limsup_{\lambda \rightarrow 0} \int_0^T \!\!\int_{t}^{T}
          ( \eps \xi_{\eps\lambda} , \dot{u}_{\eps\lambda} ) \leq\;&
                                                                     -
                                                                     \eps
                                                                     \int_0^T
                                                                     (
                                                                     \xi_{\eps}
                                                                     ,
                                                                     u_{\eps}
                                                                     )
                                                                     +
                                                                     \int_0^T
                                                                     \langle
                                                                     -\eps
                                                                     \dot{\xi}_{\eps}
                                                                     ,
                                                                     u_{\eps}
                                                                     \rangle_{H^1(
                                                                     0,
                                                                     T
                                                                     ;
                                                                    H)
                                                                     \cap
                                                                     L^2(
                                                                     0,
                                                                     T
                                                                     ;X)}.
	\end{align*}
	Finally, via the integration by parts formula of Lemma \ref{lemmaUlisse}, it follows that
	\begin{equation*}
	\limsup_{\lambda \rightarrow 0} \int_0^T\!\! \int_{t}^{T} ( \eps \xi_{\eps\lambda} , \dot{u}_{\eps\lambda} ) \leq \int_0^T \!\!\int_{t}^{T} ( \eps\xi_{\eps}, \dot{u}_{\eps} )  
	\end{equation*}
	or, equivalently,
        	\begin{equation}
        	\label{final3}
	\limsup_{\lambda \rightarrow 0} \int_0^T t( \eps
        \xi_{\eps\lambda} , \dot{u}_{\eps\lambda} ) \leq \int_0^T t ( \eps\xi_{\eps}, \dot{u}_{\eps} ) .
	\end{equation}
        
        Moving from this, we
 can now identify $\xi_{\eps}  = {\d}_2\psi(u_{\eps}
        ,\dot{u}_{\eps} )$ by following the classical argument
        from \cite[Lemma~3.57, p.~356 and Prop.~3.59, p.~361]{attouch}. Owing to 
        \eqref{7ipotesinumero} for $\mu(t)=t$, by the strong convergence
        \eqref{convergencejuc0} we can get that 
        $\int_0^T t  \psi(J_{\lambda}u_{\eps\lambda},\cdot) \rightarrow
        \int_0^T t  \psi(u_{\eps},\cdot)$ on $L^2(0,T;H)$ in the sense of
        Mosco. In particular, 
        for any $v\in L^2(0,T;H)$ there exists a recovery sequence $ v_{\lambda} \rightarrow v$ in $L^2(0,T;H)$ such that 
	\begin{equation*}
	\limsup_{\lambda \rightarrow 0} \int_0^T t 
        \psi(J_{\lambda}u_{\eps\lambda},v_{\lambda})  =   \int_0^T  t  \psi(u_{\eps},v).
      \end{equation*}
           On the other hand,   since $\dot{u}_{\eps\lambda} \rightharpoonup \dot{u}_{\eps}$ in $L^2(0,T;H)$, we get
	\begin{equation*}
	\liminf_{\lambda \rightarrow 0} \int_0^T t  \psi(J_{\lambda}u_{\eps\lambda},\dot{u}_{\eps\lambda}) \geq \int_0^T t  \psi(u_{\eps},\dot{u}_{\eps}).
	\end{equation*}
 As $\xi_{\eps\lambda} = \d_2\psi(J_\lambda u_{\eps\lambda},\dot
u_{\eps\lambda})$ one has that
\begin{equation*}\int_0^T t  \psi(J_{\lambda}u_{\eps\lambda},v_{\lambda}) -\int_0^T t  
\psi(J_{\lambda}u_{\eps\lambda},\dot{u}_{\eps\lambda}) \geq \int_0^T t
(\xi_{\eps\lambda},v_{\lambda}-
\dot{u}_{\eps\lambda}).
\end{equation*}
By passing to the $\limsup$ as $\lambda \to 0$  and by using
\eqref{final3} we obtain 
$$\int_0^Tt\psi(u_{\eps},v) -\int_0^Tt\psi(u_{\eps},\dot{u}_{\eps}) \geq \int_0^Tt(\xi_{\eps},v- \dot{u}_{\eps}).$$
 As $ v \in L^2(0,T;H)$ is arbitrary, this  corresponds to
	\begin{equation*}
	\xi_{\eps}  = {\d}_2\psi(u_{\eps} ,\dot{u}_{\eps} ) \;\;
	\text{in} \;H, \; \text{a.e. in}\ (0,T).
      \end{equation*}

$\bullet$	 {\it Step 6: Strong convergence of $\sqrt{t}\dot
  u_{\eps\lambda}$.} As 
$\xi_{\eps\lambda} = \d_2 \psi(J_\lambda
u_{\eps\lambda},\dot{u}_{\eps\lambda})$ and $\xi_\eps = \d_2 \psi(u_\eps,\dot{u}_\eps)$, 
from \eqref{str-mono}, \eqref{convergencejuc0},  and \eqref{final3} we deduce that
\begin{align*}
c_{12,R} \int^T_0 t \|\dot{u}_{\eps\lambda}-\dot{u}_\eps\|_H^2
&\leq \int^T_0 t \left(\xi_{\eps\lambda} - \xi_\eps, \dot{u}_{\eps\lambda} - \dot{u}_\eps \right)\\
&\quad + c_{13,R} \sup_{t \in [0,T]} \|J_\lambda u_{\eps\lambda}-u_\eps\|_H^2 \int^T_0 t\|\dot{u}_\eps\|_H^2 \to 0,
\end{align*}
which implies
\begin{equation}\label{eq:dotu0}
\sqrt{t} \, \dot{u}_{\eps\lambda} \to \sqrt{t} \, \dot{u}_\eps \quad \mbox{ strongly in } L^2(0,T;H).
\end{equation}
In particular, we have that
\begin{equation}\label{eq:dotu}\dot{u}_{\eps\lambda}(t) \rightarrow
\dot{u}_{\eps}(t) \quad \text{in} \ \ H, \ \text{for a.e.} \ t \in (0,T).
\end{equation}

$\bullet$	{\it Step 7: Identification of $\gamma_{\eps}$.}
As
$\dot{u}_{\eps\lambda}(t) \rightarrow
\dot{u}_{\eps}(t)$ in $H$ for almost every $t \in (0,T)$ from
\eqref{eq:dotu},
$J_{\lambda}u_{\eps\lambda}(t) \rightarrow u_{\eps}(t)$ in $H$ for
every $t\in [0,T]$ from \eqref{convergencejuc0}, and $\d_1\psi \in C^1(H\times H;H)$ we deduce that
\begin{equation*}
{\d}_1\psi(J_{\lambda}u_{\eps\lambda}(t)  ,
\dot{u}_{\eps\lambda}(t) ) \rightarrow {\d}_1\psi
(u_{\eps} (t) , \dot{u}_{\eps} (t)) \;\; \text{in} \; H, \;
\text{for a.e. $t$ in}\ (0,T).  
\end{equation*}
Eventually, since by \eqref{0ipotesinumero}
\begin{equation*}
\lvert {\d}_1 \psi(J_{\lambda}u_{\eps\lambda} , \dot{u}_{\eps\lambda})
\rvert_{H} \leq c_{5,R}(  1+\lvert \dot{u}_{\eps\lambda} \rvert_H^2   ),
\end{equation*}
we apply a generalization of the Dominated Convergence Theorem (see Lemma \ref{GDCT}) and we infer that ${\d}_1 \psi(J_{\lambda}u_{\eps\lambda} , \dot{u}_{\eps\lambda})  \rightarrow  {\d}_1 \psi(u_{\eps} , \dot{u}_{\eps}) $ in $L^1(0,T;H)$, namely
\begin{equation*}
\gamma_{\eps}  = {\d}_1\psi(u_{\eps},\dot{u}_{\eps} ) \;\;
\text{in} \; H, \; \text{a.e. in}  \ (0,T).
\end{equation*}

As $\gamma_{\eps}\in L^1(0,T;H)$, we find that $\dot{\xi}_{\eps} \in
L^1(0,T;H) + L^2(0,T;X^*)$, so that equation
\eqref{limitequationepsilon} holds in $L^1(0,T;H) + L^2(0,T;X^*)$. In conclusion, we have
\begin{equation*}
-\eps \dot{\xi}_{\eps}  + \xi_{\eps}  + \eps
\gamma_{\eps}  + \eta^1_{\eps}  +\eta^2_{\eps} 
= 0 \;\; \text{in} \; X^*,\;\text{a.e. in}\   (0,T).
\end{equation*}

\subsection{Minimization of the WED Functional $W^{\eps}$}
\label{minimization}
In this section, we show that   $u_{\eps}$   minimizes the WED
functional $W^{\eps} $.

For any given $v \in K(u_0)$, consider the regularized WED functional
\begin{equation*}
W^{\eps\lambda}(v) = \int_{0}^T e^{-t/\eps} \biggl( \eps \psi(J_{\lambda}v, \dot{v}) + \phi^1_{\lambda}(v) + \phi_{\lambda}^2(v)   \biggr).
\end{equation*}
From the basic properties of the resolvent $J_{\lambda}$ and of the
Moreau--Yosida regularizations $\phi^1_{\lambda}$ and
$\phi_{\lambda}^2$, since $J_{\lambda}v \rightarrow v$ in $H$ almost
everywhere in $ (0,T)$, by the continuity of $\psi$ we deduce
that $\psi(J_{\lambda}v,\dot{v}) \rightarrow
\psi(v,\dot{v})$ almost everywhere in $ (0,T)$, that $\phi_{\lambda}^2(v) \rightarrow \phi^2(v)$ a.e.~in $t \in (0,T)$, and that $\phi^1_{\lambda}(v) \rightarrow \phi^1(v)$ a.e.~in $t\in(0,T)$. Moreover, by \eqref{2.7ipotesinumero} it holds that
\begin{equation*}
0 \leq \eps \psi(J_{\lambda}v,\dot{v}) + \phi^1_{\lambda}(v)
+\phi_{\lambda}^2(v) \leq \eps c \left( 1+\lvert \dot{v} \rvert_H^2   \right) + \phi^1(v) + \phi^2(v) \in L^1(0,T).
\end{equation*}
Thus, we may apply the Dominated Convergence Theorem and obtain
\begin{equation}
\label{no1}
W^{\eps\lambda}(v) \rightarrow W^{\eps}(v) \;\; \text{as}\;\lambda\rightarrow 0.
\end{equation}

As $u_{\eps\lambda} \in K_0(u_0)$  is a global minimizer of $W^{\eps\lambda}$, we have
\begin{equation}
\label{no2}
	W^{\eps\lambda} (v) \geq W^{\eps\lambda}(u_{\eps\lambda}) \;\;\text{for any}\;v\in K(u_0).
\end{equation}
Convergences \eqref{convergencelambdau}, \eqref{convergenzaTemp}, and \eqref{convergencejuc0}, the lower semicontinuity of $\phi^1$ and $\phi^2$, and the convergence \eqref{7ipotesinumero} imply that
\begin{align}
&\liminf_{\lambda\rightarrow 0} W^{\eps\lambda}(u_{\eps\lambda}) = \liminf_{\lambda\rightarrow 0} \int_{0}^T e^{-t/\eps} \biggl( \eps \psi(J_{\lambda}u_{\eps\lambda},\dot{u}_{\eps\lambda}) + \phi^1_{\lambda}(u_{\eps\lambda}) +\phi_{\lambda}^2(u_{\eps\lambda})  \biggr)\nonumber\\[2mm]
&\qquad \geq \liminf_{\lambda\rightarrow 0} \int_{0}^T e^{-t/\eps} \biggl( \eps \psi(J_{\lambda}u_{\eps\lambda},\dot{u}_{\eps\lambda}) + \phi^1(J_{\lambda}u_{\eps\lambda}) + \phi^2 (I_{\lambda}u_{\eps\lambda})   \biggr)
\geq W^{\eps}(u_{\eps}).\label{no3}
\end{align}
 By combining  convergence \eqref{no1} with inequalities \eqref{no2} and \eqref{no3} we obtain
\begin{equation*}
W^{\eps}(v)=  \lim_{\lambda\rightarrow 0}  W^{\eps\lambda}(v)
\geq \liminf_{\lambda\rightarrow 0}
W^{\eps\lambda}(u_{\eps\lambda})\geq W^{\eps}(u_{\eps})\quad \text{for any}\; v\in K(u_0),
\end{equation*}
which proves that $u_{\eps} \in K(u_0)$   minimizes $W^{\eps}$.
  In case $u_\eps$ is the unique global minimizer of 
 $W^{\eps}$, we have proved that $u_{\eps}$ solves the Euler--Lagrange problem \eqref{eulerolagrangeepsilon}--\eqref{eulerolagrangeepsilon4} and Theorem \ref{theorem}.${\rm (ii)}$ follows.

In case $W^{\eps}$ has more  global  minimizers, the
convergence analysis from Section \ref{passagelimit} applies just to
those resulting as limits $\lim_{\lambda \rightarrow 0}
u_{\eps\lambda}$ of (subsequence of)  global  minimizers
$u_{\eps\lambda}$ of the regularized WED functional
$W^{\eps\lambda}$. In this case, we are forced to refine the
convergence argument by arguing by penalization. Let
$\hat{u}_{\eps} \in K(u_0)$ be a fixed  global 
minimizer of $W^{\eps}$  and define the penalized functionals $\hat{W}^{\eps}$ and $\hat{W}^{\eps\lambda}$ by
\begin{equation}
  \hat{W}^{\eps}(u) =
\left\{
  \begin{array}{ll} 
  \displaystyle \int_{0}^T e^{-t/\eps} \bigl( \eps\psi(u,\dot{u}) +
    \hat{\phi}(u)  \bigr)\,\d t \quad&\text{if}\ u \in K(u_0),\\[5mm]
                                     \infty &\text{else}
  \end{array}
\right.
\label{penalized}
\end{equation}
with
\begin{equation*}
\hat{\phi}(u) = \phi^1(u) + \phi^2(u) + \frac{1}{2} \lvert u - \hat{u}_{\eps} \rvert_H^2,
\end{equation*}
and
\begin{equation*}
\hat{W}^{\eps\lambda}(u)=\int_{0}^T e^{-t/\eps} \bigl( \eps\psi(J_{\lambda}u,\dot{u}) + \hat{\phi}_{\lambda}(u)  \bigr),
\end{equation*}
with
\begin{equation*}
\hat{\phi}_{\lambda}(u) = \phi^1_{\lambda}(u) + \phi_{\lambda}^2(u) + \frac{1}{2} \lvert u - \hat{u}_{\eps} \rvert_H^2.
\end{equation*}
From the definition \eqref{penalized} of $\hat{W}^{\eps}$ we readily check that $\hat{u}_{\eps}$ is the unique global minimizer of $\hat{W}^{\eps}$ in $K(u_0)$, and that $\hat{W}^{\eps}(\hat{u}_{\eps}) = W^{\eps}(\hat{u}_{\eps}) = \min_{v \in K(u_0)} W^{\eps}(v)$.
Arguing as in Section \ref{directmethod}, we can show that
$\hat{W}^{\eps\lambda}$ admits a  global  minimizer
$\tilde{u}_{\eps\lambda}$. Moreover, $\tilde{u}_{\eps\lambda}$
satisfies a penalized version of the regularized Euler--Lagrange
problem \eqref{euleroepsilonlambda}--\eqref{datoiniziale}, namely,
\begin{align}
&-\eps\dot{\tilde{\xi}}_{\eps\lambda}  + \tilde{\xi}_{\eps\lambda} 
  + \eps J_{\lambda}^*\tilde{\gamma}_{\eps\lambda} 
  +\tilde{\eta}^1_{\eps\lambda}  +\tilde{\eta}^2_{\eps\lambda}  =
                                                                     -
                                                                     \tilde{u}_{\eps\lambda} 
                                                                     +
                                                                     \hat{u}_{\eps} 
  \quad\text{in}\  H,  \ \text{a.e. in} \ (0,T),  \label{lagrangelambdapenalized}\\ 
&\tilde{\xi}_{\eps\lambda}  = {\d}_{2}  \psi(J_{\lambda}\tilde{u}_{\eps\lambda}   ,\dot{\tilde{u}}_{\eps\lambda}  )\quad\text{in}\  H,  \ \text{a.e. in} \ (0,T),  \\
&\tilde{\gamma}_{\eps\lambda}  = {\d}_{1}  \psi(J_{\lambda}\tilde{u}_{\eps\lambda}                    ,
                                                                                   \dot{\tilde{u}}_{\eps\lambda}                          )\quad\text{in}\
                                                                                   H,  \ \text{a.e. in} \ (0,T),  \\
&\tilde{\eta}^1_{\eps\lambda}   = A_\lambda \tilde{u}_{\eps\lambda}\quad\text{in}\ H, \ \text{a.e. in} \ (0,T),\\
&\tilde{\eta}^2_{\eps\lambda}   = \partial      \phi_{\lambda}^2(\tilde{u}_{\eps\lambda}                          )\quad\text{in}\     H,  \ \text{a.e. in} \ (0,T),  \\
&\tilde{u}_{\eps\lambda}(0)  = u_0 ,\\
& \epsi  \tilde{\xi}_{\eps\lambda}(T)  = 0 . 
\end{align} 
We now aim at recovering estimates \eqref{154}--\eqref{201} in this penalized setting. We consider the time integral over $(0,T)$ of the test of equation \eqref{lagrangelambdapenalized} on $\dot{\tilde{u}}_{\eps\lambda}$ getting
\begin{align}
\int_{0}^T \biggl( -\eps ( \dot{\tilde{\xi}}_{\eps\lambda} ,
  \dot{\tilde{u}}_{\eps\lambda})  + ( \tilde{\xi}_{\eps\lambda} ,&
                                                                   \dot{\tilde{u}}_{\eps\lambda}) + \eps ( \tilde{\gamma}_{\eps\lambda} , J_{\lambda}\dot{\tilde{u}}_{\eps\lambda}) + ( \tilde{\eta}^1_{\eps\lambda} , \dot{\tilde{u}}_{\eps\lambda} ) + ( \tilde{\eta}^2_{\eps\lambda} , \dot{\tilde{u}}_{\eps\lambda} ) \biggr)\nonumber \\&
= -  \int_{0}^T  (\tilde{u}_{\eps\lambda} - \hat{u}_{\eps} , \dot{\tilde{u}}_{\eps\lambda}) \label{reallytmep}.
\end{align}
The right-hand side term can be treated as follows\/{\rm :}
\begin{equation*}
 -\int_{0}^T  (\tilde{u}_{\eps\lambda} - \hat{u}_{\eps} , \dot{\tilde{u}}_{\eps\lambda}) = - \frac{1}{2} \int_{0}^T \frac{\d}{\d t} \lvert \tilde{u}_{\eps\lambda} \rvert_H^2 + \int_0^T (\hat{u}_{\eps},\dot{\tilde{u}}_{\eps\lambda}) .
\end{equation*}
Inserting this equality into \eqref{reallytmep} we get
\begin{align*}
\int_{0}^T \biggl( -\eps ( \dot{\tilde{\xi}}_{\eps\lambda} , \dot{\tilde{u}}_{\eps\lambda} ) + ( \tilde{\xi}_{\eps\lambda} , &\dot{\tilde{u}}_{\eps\lambda}) + \eps ( \tilde{\gamma}_{\eps\lambda} , J_{\lambda}\dot{\tilde{u}}_{\eps\lambda}) + ( \tilde{\eta}^1_{\eps\lambda} , \dot{\tilde{u}}_{\eps\lambda}) + ( \tilde{\eta}^2_{\eps\lambda} , \dot{\tilde{u}}_{\eps\lambda}) \biggr) \\
& =- \frac{1}{2} \int_{0}^T \frac{\d}{\d t} \lvert \tilde{u}_{\eps\lambda} \rvert_H^2 + \int_0^T (\hat{u}_{\eps},\dot{\tilde{u}}_{\eps\lambda}) ,
\end{align*}
which by Lemma \ref{lemmachainrule} and $\eps \tilde \xi_{\eps \lambda} (T)=0$ can be rearranged as follows:
\begin{equation*}
	\begin{aligned}
		  &\int_0^T  \frac{\d}{\d t}\biggl( \eps\psi(J_{\lambda}\tilde{u}_{\eps\lambda} , \dot{\tilde{u}}_{\eps\lambda}) + \phi^1_{\lambda}(\tilde{u}_{\eps\lambda}) + \phi_{\lambda}^2(\tilde{u}_{\eps\lambda}) + \frac{1}{2} \lvert \tilde{u}_{\eps\lambda}  \rvert_H^2 \biggr) +\int_0^T ( \tilde{\xi}_{\eps\lambda} , \dot{\tilde{u}}_{\eps\lambda} )\\
		  &\qquad+\eps ( \tilde{\xi}_{\eps\lambda}(0) , \dot{\tilde{u}}_{\eps\lambda}(0) ) =  \int_{0}^T (\hat{u}_{\eps} , \dot{\tilde{u}}_{\eps\lambda} ).
	\end{aligned}
\end{equation*}
The latter entails that
\begin{align}
&\int_0^T ( \tilde{\xi}_{\eps\lambda},\dot{\tilde{u}}_{\eps\lambda}) + \eps \biggl( ( \tilde{\xi}_{\eps\lambda}(0) , \dot{\tilde{u}}_{\eps\lambda}(0) ) - \psi(J_{\lambda}\tilde{u}_{\eps\lambda}(0) , \dot{\tilde{u}}_{\eps\lambda}(0))	\biggr)\nonumber \\
&\quad + \phi^1_{\lambda}(\tilde{u}_{\eps\lambda}(T))+ \phi_{\lambda}^2(\tilde{u}_{\eps\lambda}(T)) + \frac{1}{2} \lvert \tilde{u}_{\eps\lambda}(T)\rvert_H^2    \nonumber \\ &\leq \phi^1_{\lambda}(\tilde{u}_{\eps\lambda}(0)) + \phi_{\lambda}^2(\tilde{u}_{\eps\lambda}(0)) + \frac{1}{2} \lvert \tilde{u}_{\eps\lambda}(0)\rvert_H^2 +\delta \int_0^T \lvert \dot{\tilde{u}}_{\eps\lambda} \rvert_H^2 +\frac{1}{4\delta}\int_0^T \lvert {\hat{u}}_{\eps} \rvert_H^2,\label{temppenalized}
\end{align}
where $\delta>0$ will be chosen small, see below. From the boundedness of the terms $\phi^1_{\lambda}(\tilde{u}_{\eps\lambda}(0)) = \phi^1_{\lambda}(u_0)$, $\phi^2(\tilde{u}_{\eps\lambda}(0)) = \phi_{\lambda}^2(u_0)$, and $\lVert {\hat{u}}_{\eps} \rVert_{L^2(0,T;H)}$ and from the nonnegativity of all pointwise terms in the left-hand side of \eqref{temppenalized} we infer that
\begin{align*}
	\int_0^T ( \tilde{\xi}_{\eps\lambda},\dot{\tilde{u}}_{\eps\lambda}) \leq c + \delta \int_0^T \lvert \dot{\tilde{u}}_{\eps\lambda} \rvert_H^2.
\end{align*}
Hence, by choosing $\delta$ small enough, property
\eqref{2.7ipotesinumero} implies that 
\begin{equation}
\label{boundlambdautemp}
\int_0^T \lvert \dot{\tilde{u}}_{\eps\lambda} \rvert_H^2 \leq c.
\end{equation}
By means of the bound \eqref{boundlambdautemp}, inequality \eqref{temppenalized} leads to estimates analogous to \eqref{155}--\eqref{201}, where $(u_{\eps\lambda},\xi_{\eps\lambda},\gamma_{\eps\lambda},\eta^1_{\eps\lambda},\eta^2_{\eps\lambda})$ are replaced by their penalized counterparts \\ $(\tilde{u}_{\eps\lambda},\tilde{\xi}_{\eps\lambda},\tilde{\gamma}_{\eps\lambda},\tilde{\eta}^1_{\eps\lambda},\tilde{\eta}^2_{\eps\lambda})$.
Hence, we may extract convergent subsequences and proceed precisely as in Section \ref{passagelimit}. Eventually, we obtain, for some not relabeled subsequences, the following convergences\/{\rm :}
\begin{alignat}{4}
\tilde{u}_{\eps\lambda} &\rightharpoonup \tilde{u}_{\eps}\quad &&\text{in} \quad H^1(0,T;H),\label{0convergencelambdautemp}\\
\tilde{u}_{\eps\lambda} &\rightarrow \tilde{u}_{\eps} \quad &&\text{in}
                                                      \quad
                                                      L^{q}(0,T;H)
                                                      \quad \forall q < \infty,\label{0strongconvergencelambdautemp}\\
J_{\lambda}^*\tilde{\gamma}_{\eps\lambda} &\rightharpoonup \tilde{\gamma}_{\eps} \quad &&\text{in} \quad (H^1(0,T;H))^*\label{0convergencelambdagammatemp},\\
\tilde{\eta}^1_{\eps\lambda} &\rightharpoonup \tilde{\eta}^1_{\eps} \quad &&\text{in} \quad L^2(0,T;X^*),\label{0convergencelambdaetatemp}\\
\tilde{\eta}^2_{\eps\lambda} &\rightharpoonup \tilde{\eta}^2_{\eps}
\quad &&\text{in} \quad L^r(0,T;H),\label{0convergence}\\
\tilde{\xi}_{\eps\lambda} &\weakstar \tilde{\xi}_{\eps} \quad
                                                        &&\text{in}
                                                        \quad
                                                        {L^2(0,T;H)\cap BV([0,T];X^*)},\label{0convergencelambdaxitemp}\\
\dot{\tilde{\xi}}_{\eps\lambda} &\rightharpoonup \dot{\tilde{\xi}}_{\eps} \quad &&\text{in} \quad (H^1(0,T;H))^* + L^2(0,T;X^*),\label{0convergencelambdaxipuntotemp}
\end{alignat}
where $(\tilde{u}_{\eps},\tilde{\xi}_{\eps},\tilde{\gamma}_{\eps},\tilde{\eta}^1_{\eps},\tilde{\eta}^2_{\eps})$ solve
\begin{equation}
	-\eps\dot{\tilde{\xi}}_{\eps} + \tilde{\xi}_{\eps} +\eps\tilde{\gamma}_{\eps} +\tilde{\eta}^1_{\eps}  +\tilde{\eta}^2_{\eps}= - \tilde{u}_{\eps} + \hat{u}_{\eps}\quad \text{in} \quad (H^1(0,T;H))^* + L^2(0,T;X^*) . \label{0eulerolagrangeepsilontemp}
\end{equation}
To identify the limit functions we proceed as in Section \ref{passagelimit}, the only difference being in equation \eqref{l1temp}, namely in the identification of $\tilde{\xi}_{\eps} = \lim_{\lambda \rightarrow 0} \tilde{\xi}_{\eps\lambda}$. In this case, an additional term arises from the penalization. Indeed, we have 
\begin{align*}
&\int_{t_1}^{t_2} ( \eps \tilde{\xi}_{\eps\lambda} ,
  \dot{\tilde{u}}_{\eps\lambda} ) = \eps ( \tilde{\xi}_{\eps\lambda}(t_2) , \tilde{u}_{\eps\lambda}(t_2) ) - \eps ( \tilde{\xi}_{\eps\lambda}(t_1) , \tilde{u}_{\eps\lambda}(t_1) ) + \int_{t_1}^{t_2} ( -\eps\dot{\tilde{\xi}}_{\eps\lambda}, \tilde{u}_{\eps\lambda} )\\[2mm]
&\qquad = \eps ( \tilde{\xi}_{\eps\lambda}(t_2) , \tilde{u}_{\eps\lambda} ) - \eps ( \tilde{\xi}_{\eps\lambda}(t_1) , \tilde{u}_{\eps\lambda}(t_1) ) + \int_{t_1}^{t_2}  ( -\tilde{\xi}_{\eps\lambda}  -\eps J_{\lambda}^* \tilde{\gamma}_{\eps\lambda} - \tilde{\eta}^1_{\eps\lambda} - \tilde{\eta}^2_{\eps\lambda} -  \tilde{u}_{\eps\lambda} + \hat{u}_{\eps} , \tilde{u}_{\eps\lambda} ),
\end{align*}
where in the last equality we have used the penalized version of the regularized Euler--Lagrange equation \eqref{lagrangelambdapenalized}. To treat the additional term we note that the strong convergence \eqref{0strongconvergencelambdautemp} of $\tilde{u}_{\eps\lambda}$ to $\tilde{u}_{\eps}$ in $L^q(0,T;H)$ for all $q < \infty$ implies the strong convergence of $\tilde{u}_{\eps\lambda} - \hat{u}_{\eps}$ to $\tilde{u}_{\eps} - \hat{u}_{\eps}$ in $L^2(0,T;H)$. Hence, we have that
\begin{equation*}
\int_{t_1}^{t_2} (\tilde{u}_{\eps\lambda} - \hat{u}_{\eps}, \tilde{u}_{\eps\lambda} ) \rightarrow \int_{t_1}^{t_2}  (\tilde{u}_{\eps} - \hat{u}_{\eps} , \tilde{u}_{\eps} ).
\end{equation*}
We can hence proceed as in Section \ref{passagelimit} and 
obtain,  
\begin{align}
&\tilde{\xi}_{\eps}  =
                        {\d}_{2}\psi(\tilde{u}_{\eps}
                      ,\dot{\tilde{u}}_{\eps} )
                        \quad \text{in}\  H,  \ \text{a.e. in} \
                (0,T),  \\
&\tilde{\gamma}_{\eps}  =
                           {\d}_{1}\psi(\tilde{u}_{\eps}
                         ,\dot{\tilde{u}}_{\eps} )
                           \quad\text{in}\  H,  \ \text{a.e. in} \
                (0,T),  \\
&\tilde{\eta}^1_{\eps}  = A\tilde{u}_{\eps}   \quad \text{in} \  X^*,  \ \text{a.e. in} \
                (0,T),  \\
&\tilde{\eta}^2_{\eps}  \in \partial \phi^2(\tilde{u}_{\eps} ) \quad
                           \text{in} \  H,  \ \text{a.e. in} \
                (0,T),  \\
&\tilde{u}_{\eps}(0) =u_0, \\
& \eps  \tilde{\xi}_{\eps}(T) = 0, \label{00provatemp}
\end{align}
 as well as 
\begin{equation*}
-\eps\dot{\tilde{\xi}}_{\eps}  + \tilde{\xi}_{\eps}
+\eps\tilde{\gamma}_{\eps}  +\tilde{\eta}^1_{\eps} 
+\tilde{\eta}^2_{\eps} = - \tilde{u}_{\eps}  +
\hat{u}_{\eps} \quad \text{in} \  X^* \ \text{a.e. in} \
                (0,T) .
\end{equation*}
Eventually, $\tilde{u}_{\eps}$ minimizes $\tilde{W}^{\eps}$. Since
$\tilde{W}^{\eps}$ admits the unique global minimizer $\hat{u}_{\eps}$
we have that $\tilde{u}_{\eps} = \hat{u}_{\eps}$. By substituting
$\tilde{u}_{\eps} = \hat{u}_{\eps}$ in
\eqref{0eulerolagrangeepsilontemp}--\eqref{00provatemp} we find that
the global minimizer $\hat{u}_{\eps}$ solves the Euler--Lagrange problem \eqref{eulerolagrangeepsilon}--\eqref{eulerolagrangeepsilon4}. The proof of Theorem \ref{theorem}.${\rm (ii)}$ is hence complete.



\section{The causal limit $\eps \rightarrow 0$}
\label{causallimit}
This section is devoted to a proof of Theorem \ref{theorem}.${\rm (iii)}$. First, we provide {\it a priori} estimates on  global  minimizers $(u_{\eps})_{\eps}$ which are independent of $\eps$. Then, we consider the limit $\eps \rightarrow 0$ and show that $u = \lim_{\eps \rightarrow 0} u_{\eps}$ solves the gradient flow \eqref{doublyNonLinearCauchy1}--\eqref{doublyNonLinearCauchy4}.

As the estimates in Subsection \ref{apriorilambda} do not depend on
$\eps$, we may pass them to the $\liminf$ for $\eps \rightarrow 0$. Hence, by using the lower semicontinuity of the norm and of $\phi^1$ and $\phi^2$, it is straightforward to deduce that
\begin{align}
	 \lVert &u_{\eps} \rVert_{H^1(0,T;H) \cap L^2(0,T;X)} + \lVert \xi_{\eps} \rVert_{L^{2}(0,T;H)} + \lVert  \gamma_{\eps} \rVert_{L^1(0,T;H)} \nonumber\\ &+ \lVert Au_{\eps} \rVert_{L^2(0,T;X^*)} + \lVert \eta^2_{\eps} \rVert_{L^r(0,T;H)}+ \lVert \eps \dot{\xi}_{\eps} \rVert_{L^1(0,T;H) + L^2(0,T;X^*)} \leq c.\label{espilonboundux}
\end{align}

We aim at proving that the solutions $(u_{\eps},\xi_{\eps},\eta^1_\eps,\eta^2_{\eps})$ of the Euler--Lagrange problem \eqref{eulerolagrangeepsilon}--\eqref{eulerolagrangeepsilon4} converge to a solution $(u,\xi,Au,\eta^2)$ of the gradient flow \eqref{doublyNonLinearCauchy1}--\eqref{doublyNonLinearCauchy4}. Estimate \eqref{espilonboundux} implies that, up to not relabeled subsequences, the following convergences hold\/{\rm :}
\begin{alignat}{4}
u_{\eps} &\rightharpoonup u \quad &&\text{in} \quad H^1(0,T;H)\cap L^2(0,T;X) \label{convergenceue0T},\\
\xi_{\eps} &\rightharpoonup \xi \quad &&\text{in} \quad L^{2}(0,T;H),\label{non}\\
\eps\gamma_{\eps} &\to 0 \quad &&\text{in} \quad L^1(0,T;H),\\
Au_{\eps} &\rightharpoonup Au \quad &&\text{in} \quad L^2(0,T;X^*),\\
\eta^2_{\eps} &\rightharpoonup \eta^2 \quad &&\text{in} \quad L^{r}(0,T;H),\label{5.5bTemp}\\
\eps\dot{\xi}_{\eps} &\rightharpoonup 0 \quad &&\text{in} \quad (H^1(0,T;H))^* + L^2(0,T;X^*).\label{convboh}
\end{alignat}
As the sequence $u_{\eps}$ is bounded in $H^1(0,T;H) \cap L^2(0,T;X)$,
the  Aubin--Lions--Simon Lemma \ref{lemma:simon}  implies that
\begin{equation*}
\label{convergencejuc0espilon}
u_{\eps} \rightarrow u \;\; \text{in} \; C([0,T];H).
\end{equation*}
The above convergences allow us to pass to the limit in the Euler--Lagrange equation \eqref{eulerolagrangeepsilon} and to obtain
\begin{equation}
\label{equationboh}
\xi  + Au  + \eta^2  = 0 \quad \text{in} \ X^*,  \ \text{a.e. in} \  (0,T).
\end{equation} 

Let us now proceed with the identification of the   limits $\eta$
and $\xi$. Note preliminarily that a comparison in  
equation \eqref{equationboh} 
guarantees  that $Au  \in H$ almost everywhere in  $(0,T)$.  In
particular,  
we have that
\begin{equation*}
	Au  = \partial\phi^1(u ) \quad \text{in} \ H, \ \text{a.e.}\  (0,T),
\end{equation*}
so that equation \eqref{equationboh} actually holds in $H$.

As $u_{\eps}$ strongly converges to $u$ in $L^2(0,T;H)$ and
$\eta_{\eps}^2$ weakly converges to $\eta^2$ in $L^r(0,T;H)$, 
by applying \cite[Prop. 2.5, p. 27]{Brezis73} we find that
\begin{equation*}
	\eta  \in \partial\phi(u ) \;\; \text{in}\; H,\;\text{a.e. in}
        \ (0,T).
\end{equation*}

We now pass to the limit as $\lambda \to 0$ in inequality 
\eqref{5.23}. Recalling the convergences  
\eqref{convergencelambdaxi}, \eqref{convergencejuc0}, and
\eqref{eq:dotu0}, as well as the fact that $(\xi_{\eps\lambda},\dot
u_{\eps\lambda}) \geq 0$ a.e., for all $\delta>0$ one has that 
\begin{align*}
 & \int_\delta^T (\xi_{\eps},\dot
u_{\eps}) = \lim_{\lambda \to 0}\int_\delta^T (\xi_{\eps\lambda},\dot
u_{\eps\lambda}) \leq \lim_{\lambda \to 0}\int_0^T (\xi_{\eps\lambda},\dot
  u_{\eps\lambda})  \\
  &\quad\leq \limsup_{\lambda \to 0} \Big(-\phi^1_\lambda
  (u_{\eps\lambda}(T)) -\phi^2_\lambda (u_{\eps\lambda}(T)) +
  \phi^1_\lambda(u_0) + \phi^2_\lambda(u_0)\Big) \leq -\phi(u_\eps(T))
  + \phi(u_0).
\end{align*}
By taking $\delta \to 0$ we hence deduce that
\begin{equation*}
\int_0^T ( \xi_{\eps} , \dot{u}_{\eps} ) \leq -\phi(u_{\eps}(T)) + \phi(u_0).
\end{equation*}
Passing to the $\limsup$ in the above inequality, we obtain
\begin{align*}
	\limsup_{\eps\rightarrow 0} \int_0^T ( \xi_{\eps} , \dot{u}_{\eps} ) &\leq - \phi(u(T)) + \phi(u_0) = \int_0^T ( -Au-\eta^2 , \dot{u} ) = \int_0^T (\xi , \dot{u} ),
\end{align*}
where the last equality follows from the already established limit
equation \eqref{equationboh}. As $\xi_{\eps} \rightharpoonup \xi$ in
$L^2(0,T;H)$, $\dot{u}_{\eps} \rightharpoonup \dot{u}$ in
$L^2(0,T;H)$ and $u_{\eps} \rightarrow u$ in $C([0,T];H)$, 
\cite[Prop~2.5,
p.~27]{Brezis73} along with \eqref{7ipotesinumero} entails that
\begin{equation*}
 \xi = {\d}_2\psi(u ,\dot{u} ) \;\; \text{in} \; H,\;\text{a.e. in} \ (0,T),
\end{equation*}
and the assertion of Theorem \ref{theorem}.${\rm (iii)}$ follows.

\section{Application}\label{sec:application}

In this section, we  present an  application of our abstract theory to
a nonlocal problem. Let $T > 0$, $\Omega$ be an open $C^{1,1}$ bounded
domain in $\mathbb{R}^d$ for $d\in \Nz$, and consider
\begin{alignat}{4}
g(k \ast u) u_t  - \Delta u + \beta(u) &\ni 0 \qquad &&\text{in}\; (0,T)\times \Omega,\label{PDE1}\\
u &=0 \qquad &&\text{on}\; (0,T)\times \partial\Omega,\\
u(0,\cdot) &= u_0 \qquad &&\text{in}\;  \Omega.\label{PDE2}
\end{alignat}
 Here, the function  $g \in C^{2,1}_{\rm loc}(\mathbb{R})$ is assumed
 to be uniformly positive, namely, $ 0 < \alpha \leq g(\cdot)$ for
 some  positive constant $\alpha$. Moreover, $\beta = \partial
 \hat{\beta}$ is the subdifferential of a proper, convex lower
 semicontinuous function $\hat{\beta}: \mathbb{R} \rightarrow
 [0,\infty]$ such that $\hat{\beta}(0) =0$ with
$$
 |\beta(r)| \leq c(1 + |r|)  \ \ \text{and} \ \ r{\beta}(r)\geq \frac1c
|r|^2- c \qquad \forall r \in\mathbb{R}
$$
for some positive constant $c$. We indicate by $u_0 \in H_0^1(\Omega)$
the given initial datum. For  all $u\in L^2(\Rz^d)$ we indicate by
$\bar{u}\in L^2(\Rz^d)$ the corresponding trivial extension, namely,
$\bar{u}=u$ on $\Omega$ and $\bar{u}=0$ elsewhere. We assume to be
given a kernel 
$k \in L^2(\mathbb{R}^d)$, so that   the symbol $k \ast u$ denotes the
convolution of $k$ with  the trivial  extension $\bar{u}$, namely, $$(k \ast u)(x) =
\int_{\mathbb{R}^d} k(x - y) \bar{u}(y) \;{\rm d}y\quad \forall x \in \Omega	.$$
 Equation \eqref{PDE1}  may arise, in some specific setting, within
 the frame of the
so-called Kirchhoff-like nonlocal models, see
\cite{Caraballo,Chipot99,Chipot01,Ferreira} and references therein.

The WED functional corresponding to \eqref{PDE1}--\eqref{PDE2} is
\begin{equation*}
W^{\eps}(u) =
\begin{cases}
\displaystyle\int_0^{T}   e^{-t/\eps}   \int_{\Omega} \left( \eps g(k
  \ast u) \frac{| u_t|^2}{2} +  \frac{|\nabla u |^2}{2} +
  \hat{\beta}(u)\right)  \,{\rm d}x   \,\d t &\text{if}\ u\in K(u_0),\\[5mm]
\infty  & \text{else},
\end{cases}
\end{equation*}
 with $K(u_0)$  given by
\begin{align*}
K(u_0) &=  \{ u \in H^1(0,T;L^2(\Omega)) \cap L^2(0,T;H^1_0(\Omega)) \; : \; u(0) = u_0  \}.
\end{align*}
\begin{theorem}
	\label{theorem6.1}
	Under the above assumptions, it holds that
	\begin{enumerate}[{\rm (i)}]
		\item $W^{\eps}$ admits a  global  minimizer $u_{\eps} \in K(u_0)$
		\item $u_{\eps} \rightarrow u$ (up to subsequences) weakly in $ H^1(0,T;L^2(\Omega)) \cap L^2(0,T;H_0^1(\Omega)) $ where $u$ solves \eqref{PDE1}--\eqref{PDE2} in a strong sense, that is, $u \in L^2(0,T;H^2(\Omega))$ and there exists $\xi \in L^2(0,T;L^2(\Omega))$ such that $\xi \in \beta(u)$ a.e.~in $(0,T) \times \Omega$ and 
$$ g(k \ast u) u_t  - \Delta u + \xi = 0 \qquad \text{a.e. in}\; (0,T)\times \Omega. $$
	\end{enumerate}
\end{theorem}
\begin{proof}
	A variational formulation of \eqref{PDE1}--\eqref{PDE2} is
        obtained by letting $H=L^2(\Omega)$, $X=H_0^1(\Omega)$, and
        by  defining 
$$ \psi(u,v) = \int_{\Omega} g(k\ast u) \frac{|v|^2}{2}	\, {\rm d}x,
\quad \phi^1(u) =  \frac12\int_{\Omega}  | \nabla u |^2 \,{\rm d}x, \quad \phi^2(u) =  \int_{\Omega} \hat{\beta}(u)	\,{\rm d}x	. $$
	In order to prove the statement, it suffices to show that
        Theorem \ref{theorem} is applicable. By letting $A := \partial
        \phi_1$ we readily check that $A(u) = -\Delta u $ and
        $\beta(u) = \partial \phi^2 (u)$  and $\eta^2 \in
        \partial\phi^2(u)$ if and only if $\eta^2 \in \beta(u)$ almost
        everywhere  for $u,\, \eta^2 \in H$.  Assumption (A1)
        is  clearly  satisfied.
	
 Let us hence concentrate on  Assumption (A2). To begin with, we
check that $\psi \in C^{0}(H \times H;\mathbb{R}_+)$. Take a sequence $ ( u_n,v_n	)_n \subset H \times H$ with
$(u_n,v_n) \rightarrow (u,v)$ in $H \times H$ and use the Young
inequality to bound 
\begin{align}
  &|\psi(u_n,v_n) - \psi(u,v)| \leq   |\psi(u_n,v_n) - \psi(u_n,v)|+
                                |\psi(u_n,v) - \psi(u,v)|\nonumber\\
  &\quad \leq \sup_{|\rho| \leq \|k\|_{L^2(\Rz^d)} \|u_n\|_H} g(\rho) \left| \int_\Omega
    \left(\frac{|v_n|^2}{2}  - \frac{|v|^2}{2}\right) \d x\right|
    \nonumber\\
  &\qquad +
   \sup_{|\rho| \leq \|k\|_{L^2(\Rz^d)} \|u_n\|_H} |g'(\rho)|\| k\ast (u_n - u)\|_{H}
    \left( \int_\Omega \frac{|v|^2}{2}\, \d x\right)\nonumber\\
    &\quad \leq \sup_{|\rho| \leq \|k\|_{L^2(\Rz^d)} \|u_n\|_H} g(\rho) \,\frac12 \Big| \| v_n\|^2_H - \|
      v \|^2_H\Big| \nonumber\\
  &\qquad+
     \sup_{|\rho| \leq \|k\|_{L^2(\Rz^d)} \|u_n\|_H}  |g'(\rho)| \,\frac12 \| k\|_{L^2(\Rz^d)}\|u_n -
      u\|_{H}\|v\|^2_H \to 0.\label{argum}
\end{align}

Let us now turn to the proof that $\psi \in C^{1}(H \times
H;\mathbb{R}_+)$. We start by computing  the G\^ateaux derivative $ {\rm d}_{G,1}
        \psi(u,v) $ in $(u,v)$ of $\psi$ with respect to the first
        variable along the direction $w \in H$  getting 
	$$ ( {\rm d}_{G,1} \psi(u,v) , w ) = \int_{\Omega} g'(k \ast
        u) \frac{| v |^2}{2} k \ast w \;{\rm d}x.	 $$
 Arguing as above one readily proves that the  function $u
\mapsto ( {\rm d}_{G,1} \psi(u,v) , w ) $ is continuous in $H$. 
On the other hand, as the map  $w \mapsto ( {\rm d}_{G,1}
\psi(u,v) , w ) $ is linear and continuous,  one concludes that
the G\^ateaux derivative ${\rm d}_{G,1}\psi$ and the partial Fr\'echet
differential ${\rm d}_1 \psi$ coincide (and are continuous).  Analogously, we compute the G\^ateaux derivative of $\psi(u,v)$ with respect to the second variable along the direction $w\in H$, namely,
	$$ ( {\rm d}_{G,2} \psi(u,v) , w ) = \int_{\Omega} g(k \ast u)
        v w \;{\rm d}x.	 $$
         Both functions $v \mapsto ( {\rm d}_{G,2} \psi(u,v) , w )
        $ and $w \mapsto ( {\rm d}_{G,2} \psi(u,v) , w )
        $  are  linear and continuous in $H$.  This proves
        that $\psi$ is Fr\'echet differentiable with respect to the
        second variable with a continuous partial Fr\'echet
        differential, namely, ${\rm d}_2
        \psi =  {\rm d}_{G,2}\psi$. 

         In order to prove that  $\psi \in C^2(H\times
        H;\mathbb{R}_+)$, we compute the second G\^ateaux derivatives of
        $\psi$. Owing to the fact that $g$ is twice differentiable, by
 fixing two directions $w,\,z\in H$ we have
        \begin{align*}
          &( {\rm d}_{G,11} \psi(u,v)w , z ) = \int_{\Omega} g''(k
          \ast u) \frac{\vert v \vert^2}{2} (k \ast z) (k \ast w)
          \;{\rm d}x,\\
          &( {\rm d}_{G,21} \psi(u,v)w , z ) = \int_{\Omega} g'(k
            \ast u)  (k \ast z) v w \;{\rm d}x, \\
           &( {\rm d}_{G,12} \psi(u,v)w , z ) = \int_{\Omega}
             g'(k \ast u)   (k \ast w)vz  \;{\rm d}x,\\
          &( {\rm d}_{G,22} \psi(u,v)w, z ) = \int_{\Omega} g (k \ast u) wz \;{\rm d}x.  
        \end{align*}
Again, by suitably modifying the argument of \eqref{argum}, taking
into account that $g''\in C^{0,1}(\Rz)$ and $k \in L^2(\Rz^d)$ one
can prove that $(u,v) \mapsto ( {\rm d}_{G,ij} \psi(u,v)w , z ) $ is
continuous in $H\times H$ for all $(w,z)\in H\times H$ and all
$i,\,j=1,2$. At the same time, $(w,z)
\mapsto ( {\rm d}_{G,ij} \psi(u,v)w , z ) $ is linear and
continuous for all $(u,v)\in H\times H$ and $i,\,j=1,2$. Hence, $\psi$ is twice Fr\'echet differentiable, with
continuous second Fr\'echet differential $\d^2 \psi$.

For all $u\in H$, the map $v \in H \mapsto \int_{\Omega} g(k\ast u)
|v|^2/2	\, {\rm d}x $ is convex and we have $\psi(u,0) = 0$.

 We now proceed in proving conditions
\eqref{3ipotesinumero}--\eqref{22ipotesinumero}.
Condition \eqref{3ipotesinumero} follows straightforwardly from the
fact that $g(\cdot) \geq \alpha >0$ and
$$\psi(u,v) \geq \frac{\alpha}{2} \| v\|^2_H$$
so that one can choose $c_4=2/\alpha$. To check inequality
\eqref{0ipotesinumero} we compute
$$ \lVert {\rm d}_1 \psi(u,v) \rVert_H = \sup_{\lvert w \rvert_H \neq
  0}\frac{( {\rm d}_{G,1} \psi(u,v) , w )}{\lvert w \rvert_H} =
\sup_{\lvert w \rvert_H \neq 0}\frac{1}{\lvert w \rvert_{H}}
\int_{\Omega} g'(k \ast u) \frac{| v |^2}{2} k \ast w \;{\rm
  d}x.	$$
Assuming $\| u \|_H\leq R$ and using the Young inequality
we get
$$\lVert {\rm d}_1 \psi(u,v) \rVert_H \leq \sup_{\lvert w \rvert_H
  \neq 0}\frac{1}{\lvert w \rvert_H}\biggl(\sup_{|\rho|\leq \lvert k
  \rvert_{L^2(\mathbb{R}^d)} R}|g'(\rho)|\biggr) \frac{1}{2} \lvert v
\rvert_H^2 \lvert k \rvert_{L^2(\mathbb{R}^d)} \lvert w \rvert_H. $$
Condition \eqref{0ipotesinumero} follows then by choosing $c_{5,R}:=
\frac{1}{2} \lvert k \rvert_{L^2(\mathbb{R}^d)} \sup_{|\rho|\leq
  \lvert k \rvert_{L^2(\mathbb{R}^d)} R}|g'(\rho)|	$. Similarly
we have  that
	\begin{align*}
&	\lvert {\rm d}_1\psi(u_1,v) - {\rm d}_1 \psi(u_2,v) \rvert_H =
                         \sup_{\lvert w \rvert_H \neq 0}\frac{( {\rm
                         d}_{G,1} \psi(u_1,v) - {\rm d}_{G,1}
                         \psi(u_2,v) , w) }{\lvert w \rvert_H}
          \\
          &\quad =\sup_{\lvert w \rvert_H \neq 0}\frac{1}{\lvert w \rvert_{H}}\int_{\Omega} ( g'(k \ast u_1) - g'(k \ast u_2)) \frac{| v |^2}{2} k \ast w \;{\rm d}x.
	\end{align*}
	As $g' $ is locally Lipschitz continuous, by using again the Young inequality we get
	\begin{align*}
	&\lvert {\rm d}_1\psi(u_1,v) - {\rm d}_1 \psi(u_2,v)
   \rvert_H\\ &\quad \leq \sup_{ \lvert w \rvert_H \neq 0}\frac{1}{\lvert w
                    \rvert_{H}}\int_{\Omega} \sup_{|\rho|\leq 2 \lvert
                    k \rvert_{L^2(\mathbb{R}^d)}R} |g''(\rho)| | k
                    \ast (u_1 - u_2) | \frac{| v |^2}{2} |  k \ast w
          | \;{\rm d}x \\
          &\quad \leq  \sup_{ \lvert w \rvert_H \neq
            0}\frac{\sup_{|\rho|\leq 2 \lvert k
            \rvert_{L^2(\mathbb{R}^d)}R} |g''(\rho)|}{\lvert w
            \rvert_{H}} \lvert k \ast (u_1 - u_2)
            \rvert_{L^{\infty}(\Omega)} \frac{\lvert v \rvert_H^2}{2}
            \lvert k \ast w \rvert_{L^{\infty}(\Omega)}\\
          &\quad \leq \sup_{ \lvert w \rvert_H \neq 0}\frac{\sup_{|\rho|\leq 2 \lvert k \rvert_{L^2(\mathbb{R}^d)}R} |g''(\rho)|}{\lvert w \rvert_{H}} \lvert k \rvert_{L^2(\mathbb{R}^d)} \lvert u_1 - u_2 \rvert_H \frac{\lvert v \rvert_H^2}{2} \lvert k \rvert_{L^2(\mathbb{R}^d)} \rvert w \rvert_H .
	\end{align*}
	Then, by choosing $c_{6,R}:= \frac{1}{2}\sup_{|\rho|\leq 2
          \lvert k \rvert_{L^2(\mathbb{R}^d)}R} |g''(\rho)| \lvert k
        \rvert_{L^2(\mathbb{R}^d)}^2	$, we get
        \eqref{1ipotesinumero}.

        	Let us now consider	
	\begin{align*}
		\lvert {\rm d}_2\psi(u,v) \rvert_H  =& \sup_{ \lvert w \rvert_H \neq 0} \frac{( {\rm d}_{G,2} \psi(u,v) , w )}{\lvert w \rvert_H} =  \sup_{ \lvert w \rvert_H \neq 0} \frac{1}{\lvert w \rvert_H} \int_{\Omega} g (k \ast u) v w \; {\rm d}x.
	\end{align*}
The Young inequality and $\lvert u \rvert_H \leq R$ yield 
	\begin{align*}
		\lvert {\rm d}_2\psi(u,v) \rvert_H \leq& \sup_{ \lvert w \rvert_H \neq 0} \frac{1}{\lvert w \rvert_H} \lvert g( k \ast u) \rvert_{L^{\infty}(\Omega)} \lvert v \rvert_H \lvert w \rvert_H \leq \sup_{\rho \leq \lvert k \rvert_{L^2(\mathbb{R}^d)} R} g(\rho) \,\lvert v \rvert_H,
	\end{align*}
	so that \eqref{2ipotesinumero} follows by choosing  $c_{7,R}:=  
        \sup_{\rho \leq \lvert k \rvert_{L^2(\mathbb{R}^d)} R}
        |g(\rho)|^2		$.

       We now proceed toward \eqref{23ipotesinumero} by computing 
		\begin{align*}
		&	\lvert {\rm d}_2\psi(u_1,v) - {\rm d}_2
                  \psi(u_2,v) \rvert_H =  \sup_{\lvert w \rvert_H
                                          \neq 0}\frac{( {\rm d}_{G,2}
                                          \psi(u_1,v) - {\rm d}_{G,2}
                                          \psi(u_2,v) , w) }{\lvert w
                                          \rvert_H} \\
                  &\quad =\sup_{\lvert w \rvert_H \neq 0}\frac{1}{\lvert w \rvert_{H}}\int_{\Omega} ( g(k \ast u_1) - g(k \ast u_2)) v w \;{\rm d}x.
		\end{align*}
As $g $ is locally Lipschitz continuous, by applying
        the Young inequality we get
	\begin{align*}
&	\lvert {\rm d}_2\psi(u_1,v) - {\rm d}_2 \psi(u_2,v) \rvert_H\\
          &\quad \leq \sup_{ \lvert w \rvert_H \neq 0}\frac{1}{\lvert w \rvert_{H}}\int_{\Omega} \sup_{|\rho|\leq 2 \lvert k \rvert_{L^2(\mathbb{R}^d)}R} |g'(\rho)| | k \ast (u_1 - u_2) | | v | | w  | \;{\rm d}x \\
&\quad 	 \leq  \sup_{ \lvert w \rvert_H \neq 0}\frac{\sup_{|\rho|\leq 2 \lvert k \rvert_{L^2(\mathbb{R}^d)}R} |g'(\rho)|}{\lvert w \rvert_{H}} \lvert k \ast (u_1 - u_2) \rvert_{L^{\infty}(\Omega)} \lvert v \rvert_H\lvert w \rvert_{H}\\ 
&\quad 	\leq  \sup_{|\rho|\leq 2 \lvert k \rvert_{L^2(\mathbb{R}^d)}R} |g'(\rho)| \lvert k \rvert_{L^2(\mathbb{R}^d)} \lvert u_1 - u_2 \rvert_H  \lvert v \rvert_H,
	\end{align*}
        so that \eqref{23ipotesinumero} follows by choosing   $c_{8,R}:= \sup_{|\rho|\leq 2 \lvert k \rvert_{L^2(\mathbb{R}^d)}R} |g'(\rho)| \lvert k \rvert_{L^2(\mathbb{R}^d)}	$.
	
We now consider
		\begin{align*}
		&	\lVert \d_{21} \psi (u,v) w \rVert_H =  \sup_{ \lvert z \rvert_H \neq 0} \frac{( {\rm d}_{G,21} \psi(u,v)w , z )}{\lvert z \rvert_H}\\
			&\quad =  \sup_{ \lvert z \rvert_H \neq 0} \frac{1}{\lvert z \rvert_H} \int_{\Omega} g' (k \ast u) (k \ast z) v w \; {\rm d}x.
		\end{align*}
		Thanks to the Young inequality and to $\lVert u \rVert \leq R$ we obtain
		\begin{align*}
	&	\lvert {\rm d}_{21}\psi(u,v) w \rvert_H \leq \sup_{ \lvert z \rvert_H \neq 0} \frac{1}{\lvert z \rvert_H} \lvert g'( k \ast u) \rvert_{L^{\infty}(\Omega)} \lVert k \ast z \rVert_{L^{\infty}(\Omega)} \lvert v \rvert_H \lvert w \rvert_H \\
	&\quad 	\leq \sup_{\rho \leq \lvert k
   \rvert_{L^2(\mathbb{R}^d)} R} |g'(\rho)| \lvert k
   \rvert_{L^2(\mathbb{R}^d)}  \lvert v \rvert_H \lvert w \rvert_H,
		\end{align*}
		so that \eqref{21ipotesinumero} follows for $c_{9,R}:=	\sup_{\rho \leq \lvert k \rvert_{L^2(\mathbb{R}^d)} R} |g'(\rho)| \lVert k \rvert_{L^2(\mathbb{R}^d)}$.
By computing
		\begin{align*}
			(\d_{22} \psi (u,v)w, w) = \int_{\Omega} g(k \ast u) \vert w \vert^2 \; {\rm d}x
		\end{align*}
	  and recalling that $g(\cdot) \geq \alpha > 0$, one obtains \eqref{22ipotesinumero} with $c_{10,R}:=\alpha$.


Having checked Assumption (A2), we are in the position of applying Theorem \ref{theorem} which concludes the proof. 
\end{proof}

Before closing this section, let us remark that the abstract theory applies to other classes of PDEs, as well. For instance one could consider the PDE problem
	\begin{align*}
	g( \lvert u \rvert_{L^2(\Omega)}^2) u_t  - \Delta u + \beta(u) &\ni 0 \qquad \text{a.e. in}\; (0,T)\times \Omega,\\
	u &=0 \qquad \text{a.e. on}\; (0,T)\times \partial\Omega,\\
	u(0,x) &= u_0(x) \qquad \text{a.e. in}\;  \Omega.
	\end{align*}
	Under the same assumptions above we can prove an analogous version of Theorem \ref{theorem6.1}.

        \subsection*{\it Acknowledgement}
 This research was funded
in whole or in part by the Austrian Science Fund (FWF) grants 10.55776/F65, 10.55776/I4354, 10.55776/I5149, 10.55776/P32788, and
10.55776/W1245. For
open access purposes, the authors have applied a CC BY public copyright
license to any author accepted manuscript version arising from this
submission. US acknowledges 
financial support from  the OeAD-WTZ grant CZ  05/2024.  GA is also supported by JSPS KAKENHI Grant Numbers JP21KK0044, JP21K18581, JP20H01812 and JP20H00117 as well as the Research Institute for Mathematical Sciences, an International Joint Usage/Research Center located in Kyoto University.


\appendix

\section{Generalized Dominated Convergence Theorem}
\label{appendix}

For the reader's convenience, we record here the
generalized Dominated Convergence Theorem
from \cite[Rem.~(19a), p.~1015]{zeidler}.

\begin{lemma}
	\label{GDCT}
  Let $f_n,  g_n \in L^1(\Omega)$ with
  $|f_n|\leq  g_n$ a.e., $f_n \to f$ pointwise a.e., $g_n \to g$ in $L^1(\Omega)$. Then
  $f_n \to f$ in $L^1(\Omega)$.
\end{lemma}

\begin{proof} As $g_n \to g$ a.e., by passing to the limit in $|f_n|\leq g_n$ one  has that $|f|\leq g$ a.e., so that $f\in L^1(\Omega)$.
  Since $0\leq g+g_n - |f-f_n|$ and $\liminf_n(g+g_n - |f-f_n|)=2g$ a.e. one can apply the Fatou's Lemma and get
  $$\int_\Omega 2g \leq \liminf_n \int_\Omega (g+g_n - |f-f_n|) =
  \int_\Omega 2g - \limsup_n\int_\Omega |f-f_n|.$$
  We hence have that
  $$\limsup_n \int_\Omega |f-f_n|=0.\qedhere$$
\end{proof}



\begin{thebibliography}{99}

  
\bibitem{marveggio}
  	G.~Akagi, V.~B\"ogelein, A.~Marveggio, U.~Stefanelli. Weighted Inertia-Dissipation-Energy approach to
          doubly nonlinear wave equations. Submitted,
          2024. arXiv:2401.08856.
          
	
	\bibitem{AkagiMelchionna}
	G.~Akagi, S.~Melchionna. Elliptic-regularization of nonpotential perturbations of doubly-nonlinear gradient flows of nonconvex energies: a variational approach. {\it J. Convex Anal.} 25 (2018), no. 3, 861--898.
	
	\bibitem{AkagiMelchionnaStefanelli}
	G.~Akagi, S.~Melchionna, U.~Stefanelli. Weighted energy-dissipation approach for doubly-nonlinear problems on the half line. {\it J. Evol. Equ.} 18 (2017), 49--74.
	
	\bibitem{Akagi1}
	G.~Akagi, U.~Stefanelli. A variational principle for doubly nonlinear evolution. {\it Appl. Math. Lett.} 23 (2010), 1120--1124.
	
	\bibitem{Akagi2}
	G.~Akagi, U.~Stefanelli. Weighted energy-dissipation functionals for doubly nonlinear evolution. {\it J. Funct. Anal.} 260 (2011), 2541--2578.
	
	\bibitem{Akagi3}
	G.~Akagi, U.~Stefanelli. Doubly nonlinear evolution equations as convex minimization. {\it SIAM J. Math. Anal.} 46 (2014), 1922--1945.
	
	\bibitem{Akagi4}
	G.~Akagi, U.~Stefanelli. A variational principle for gradient
        flows of nonconvex energies. {\it J. Convex Anal.} 23 (2016),
        53--75.

        

      \bibitem{attouch}
     H.~Attouch. {\it Variational convergence for functions and
       operators}. Applicable Mathematics Series. Pitman (Advanced
     Publishing Program), Boston, MA, 1984.
     

	
	
	\bibitem{bathorystefanelli}
	M.~Bathory, U.~Stefanelli. Variational resolution of outflow
        boundary conditions for incompressible Navier-Stokes. {\it
          Nonlinearity}, 35 (2022), no. 11, 5553-–5592.

        
  
\bibitem{Boegelein-et-al14}
{V.~B\"ogelein, F.~Duzaar,   P.~Marcellini}.
{Existence of evolutionary variational solutions via the calculus of variations}.
{\it J. Differential Equations}. 256 (2014),   3912--3942.
	
	\bibitem{Brezis73}
	H.~Brezis. \newblock{\it Op\'erateur maximaux monotones et semi-groupes de contractions dans les espaces de Hilbert}. \newblock (French) North-Holland Mathematics Studies, no. 5. Notas de Matem\'atica (50). North-Holland Publishing Co., Amsterdam-London; American Elsevier Publishing Co., Inc., New York, 1973.
	

	\bibitem{Caraballo}
	T.~Caraballo, M.~Herrera-Cobos, P.~Mar\'in-Rubio. Asymptotic behaviour of nonlocal $p$-Laplacian reaction-diffusion problems. {\it J. Math. Anal. Appl.} 459 (2018), 997--1015.
	
	


	\bibitem{Chipot99}
	M.~Chipot, B.~Lovat. On the asymptotic behaviour of some nonlocal problems. {\it Positivity}, 3 (1999), 65--81.


	


		\bibitem{Chipot01}
		M.~Chipot, L.~Molinet. Asymptotic behaviour of some nonlocal diffusion problems. {\it Appl. Anal.} 80 (2001), 279--315.
		
		

	
	\bibitem{ContiOrtiz}
	S.~Conti, M.~Ortiz. Minimum principles for the trajectories of
        systems governed by rate problems. {\it
          J. Mech. Phys. Solids}, 56 (2008), 1885--1904.

        

\bibitem{DalMaso93}
{G.~{Dal Maso}}.
\newblock {\it An introduction to {$\Gamma$}-convergence}.
\newblock Progress in Nonlinear Differential Equations and their Applications,
  8. Birkh\"auser Boston Inc., Boston, MA, 1993.

	
	\bibitem{deLucaDalMaso}
	G.~Dal~Maso, L.~De~Luca. A minimization approach to the wave equation on time-dependent domains. {\it Adv. Calc. Var.} 13 (2018), 425--436.
	
	
	\bibitem{DavoliStefanelli}
	E.~Davoli, U.~Stefanelli. Dynamic perfect plasticity as convex
        minimization. {\it SIAM J. Math. Anal.} 51 (2019), no. 2,
        672--730.

        
\bibitem{DeGiorgi1996}
{E.~{De Giorgi}}. {Conjectures concerning some evolution problems}.
{\it  Duke Math. J.} 81 (1996),  255--268. 
	
	\bibitem{Evans}
	L.~C.~Evans.
	\newblock{\it Partial differential equations}.
	\newblock Graduate Studies in Mathematics, volume 19. American Mathematical Society, Providence, 1998.
	
	

	\bibitem{Ferreira}
	J.~Ferreira, H.~B.~De~Oliveira. Parabolic reaction-diffusion systems with nonlocal coupled diffusivity terms. {\it Discrete Contin. Dyn. Syst.} 37 (2017), no. 5, 2431--2453.



	
	
	\bibitem{Ilmanen}
	T.~Ilmanen.
	\newblock{\it Elliptic regularization and partial regularity for motion by mean curvature}.
	\newblock Memoirs of the American Mathematical Society, number
        520. American Mathematical Society, Providence, 1994.

        


\bibitem{Hirano94}
{N. Hirano}.
{Existence of periodic solutions for nonlinear evolution equations in
Hilbert spaces}.
{\it Proc. Amer. Math. Soc.}  120 (1994),  185--192.

	
	\bibitem{LarsenOrtizRichardson}
	C.~J.~Larsen, M.~Ortiz, C.~L.~Richardson. Fracture paths from front kinetics: relaxation and rate independence. {\it Arch. Ration. Mech. Anal.} 193 (2009), 539--583. 
	
	\bibitem{LieroMelchionna}
	M.~Liero, S.~Melchionna. The weighted energy-dissipation principle and evolutionary $\Gamma$-convergence for doubly nonlinear problems. {\it ESAIM Control Optim. Calc. Var.} 25 (2019), no. 36, 1--38.
	
	\bibitem{Liero1}
	M.~Liero, U.~Stefanelli. A new minimum principle for Lagrangian mechanics, {\it J. Nonlinear Sci.} 23 (2013), 179--204.
	
	\bibitem{Liero2}
	M.~Liero, U.~Stefanelli. Weighted inertia-dissipation-energy functionals for semilinear equations. {\it Boll. Unione Mat. Ital. (9)}, 6 (2013), 1--27.
	
	\bibitem{melchionna2}
	S.~Melchionna. A variational approach to symmetry, monotonicity, and comparison for doubly-nonlinear equations.	{\it J. Math. Anal. Appl.} 456 (2017), 1303--1328.
	
	\bibitem{Melchionna201sei}
	S.~Melchionna. A variational principle for nonpotential perturbations of gradient flows of nonconvex energies. {\it J. Differential Equations}, 262 (2017), 3737--3758.
	
	\bibitem{MielkeOrtiz}
	A.~Mielke, M.~Ortiz. A class of minimum principles for
        characterizing the trajectories and the relaxation of
        dissipative systems. {\it ESAIM Control Optim. Calc. Var.} 14
        (2008), 494--516.

      \bibitem{mr}
     	A.~Mielke, T. Roub\'\i \v cek. {\it Rate-independent systems. Theory and application}. Applied Mathematical Sciences, 193. Springer, New York, 2015.
	
	\bibitem{Mielke1}
	A.~Mielke, U.~Stefanelli. A discrete variational principle for rate-independent evolution. {\it Adv. Calc. Var.} 1 (2008), 399--431.
	
	\bibitem{Mielke2}
	A.~Mielke, U.~Stefanelli. Weighted energy-dissipation functionals for gradient flows. {\it ESAIM Control Optim. Calc. Var.} 17 (2011), 52--85.
	
	
	\bibitem{ortizschmidtstefanelli}
	M.~Ortiz, B.~Schmidt, U.~Stefanelli. A variational approach to Navier-Stokes. {\it Nonlinearity}, 31 (2018), 5664--5682.
	
	\bibitem{rossisavaresegattistefanelli1}
	R.~Rossi, G.~Savaré, A.~Segatti, U.~Stefanelli. A variational principle for gradient flows in metric spaces. {\it C. R. Math. Acad. Sci. Paris}, 349 (2011), 1224--1574.
	
	\bibitem{rossisavaresegattistefanelli2}
	R.~Rossi, G.~Savaré, A.~Segatti, U.~Stefanelli. Weighted energy-dissipation functionals for gradient flows in metric spaces. {\it J. Math. Pures Appl. (9)}, 129 (2019), 1--66.
	

	
	\bibitem{scarpastefanelli1}
	L.~Scarpa, U.~Stefanelli. Stochastic PDEs via convex minimization. {\it Comm. Partial Differential Equations}, 46 (2021), no. 1, 66--97
	
	\bibitem{scarpastefanelli2}
	L.~Scarpa, U.~Stefanelli. The energy-dissipation principle for stochastic parabolic equations. {\it Adv. Math. Sci. Appl.} 30 (2021), no. 2, 429-–452.
	
	\bibitem{serratilli2}
	E.~Serra, P.~Tilli. A minimization approach to hyperbolic Cauchy problems. {\it J. Eur. Math. Soc.} 18 (2016), 2019--2044.
	
	\bibitem{serratilli1}
	E.~Serra, P.~Tilli. Nonlinear wave equations as limits of convex minimization problems: proof of a conjecture by De Giorgi. {\it Ann. of Math. (2)}, 175 (2012), 1551--1574.

	\bibitem{simon}
	J.~Simon. Compact sets in the space $L^p(0,T;B)$. {\it Ann. Mat. Pura Appl. (4)}, 146 (1987), 65--96.





	\bibitem{spadarostefanelli}
		E.~N.~Spadaro, U.~Stefanelli. A variational view at mean curvature evolution for linear growth functionals. {\it J. Evol. Equ.} 11 (2011), 793--809.

	\bibitem{stefanelli}
	U.~Stefanelli. The De Giorgi conjecture on elliptic regularization. {\it Math. Models Methods Appl. Sci.} 21 (2011), 1377--1394.
	
	
	\bibitem{tentarellitilli3}
	L.~Tentarelli. On the extensions of the De Giorgi approach to nonlinear hyperbolic equations. {\it Rend. Semin. Mat. Univ. Politec. Torino}, 74 (2016), no. 3-4, 151--160.
	
	\bibitem{tentarellitilli}
	L.~Tentarelli, P.~Tilli. An existence result for dissipative nonhomogeneous hyperbolic equations via a minimization approach. {\it J. Differential Equations}, 266 (2019), no. 8, 5185--5208.
	
	\bibitem{tentarellitilli2}
	L.~Tentarelli, P.~Tilli. De Giorgi’s approach to hyperbolic Cauchy problems: the case of nonhomogeneous equations. {\it Comm. Partial Differential Equations}, 43 (2018), no. 4, 677--698.
	
	
	\bibitem{zeidler}
	E.~Zeidler.
	\newblock{\it Nonlinear functional analysis and its applications II/B}.
	\newblock Springer-Verlag, New York, 1990.
	
\end{thebibliography}
\end{document}